\documentclass[letter,12pt]{article}
\usepackage{amsfonts}
\usepackage{amssymb}
\usepackage{amsmath}
\usepackage{cite}
\newtheorem{teo}{Theorem}
\newtheorem{lema}[teo]{Lemma}
\newtheorem{prop}[teo]{Proposition}

\newtheorem{obs}[teo]{Remark}
\newtheorem{ejem}[teo]{Example}

\newcommand{\R}{\mathbb{R}}

\newcommand{\osc}{ \mathrm{osc}}

\newcommand{\N}{\mathbb{N}}

\newcommand{\esp}{\hspace{0.05cm}}

\newcommand{\B}{\mathrm{Ball}}

\newcommand{\clo}{\mathrm{S}^1}

\begin{document}

\title{On bounded cocycles of isometries over a minimal dynamics}

\author{\Large{Daniel Coronel, \hspace{0.02cm} 
Andr\'es Navas \hspace{0.02cm} \& \hspace{0.02cm} Mario Ponce}}
\maketitle

\vspace{-0.4cm}

\begin{small}
\noindent{\bf Abstract.} 
We show the following geometric generalization of a classical theorem of W.H.~Gottschalk 
and G.A.~Hedlund: a skew action induced by a cocycle of (affine) isometries of a Hilbert space 
over a minimal dynamics has a continuous invariant section if and only if the cocycle is 
bounded. Equivalently, the associated twisted cohomological equation has a continuous 
solution if and only if the cocycle is bounded. We interpret this as a version of the 
Bruhat-Tits center lemma in the space of continuous functions. Our result also holds 
when the fiber is a proper CAT(0) space. One of the appplications concerns matrix cocycles. 
Using the action of $\mathrm{GL} (n,\mathbb{R})$ on the (nonpositively curved) space of 
positively definite matrices, we show that every bounded linear cocycle over a minimal 
dynamics is cohomologous to a cocycle taking values in the orthogonal group.
\end{small}

\tableofcontents

\newpage

\section{Introduction}

\hspace{0.45cm} Over the last years, the study of cocycles has been a central 
subject in many branches of Mathematics including not only Dynamical Systems 
and Group Theory, 
but also Geometry, 
Foliations 
and Mathematical Physics. 
This work uses ideas and techniques coming from the former two areas to deal with cocycles above 
a minimal dynamics and taking values in the group of isometries of a nonpositively curved space.  

In a general form, a {\em cocycle} associated to a dynamics on a base space is a map 
into a group $G$ that is equivariant with respect to this dynamics. These data naturally 
induce a skew action on a (perhaps nontrivial) fiber bundle, where the fibers are isomorphic 
to the phase space of the action of $G$. The possibility of ``reducing'' this fibered dynamics 
is related to a central problem, namely solving an associated {\em cohomological equation}. 
Since we are interested in the possibility of reducing our cocycles into cocycles taking 
values in some compact group, we concentrate on skew actions satisfying a natural geometric 
counterpart, namely, a boundedness property. Before stating our main (somewhat technical) 
result, we prefer to illustrate its consequences giving several applications.

\vspace{0.5cm}

\noindent{\bf A matrix version of the Gottschalk-Hedlund theorem.} Let $\Gamma$ 
be a semigroup acting minimally by homeomorphisms of a compact metric space $X$. 
Let $A$ be a linear cocycle above this action, that is, a continuous map 
$A \!: \Gamma\times X \to \mathrm{GL} (n,\mathbb{R})$ satisfying \esp 
$A(fg,x) = A(f,g (x)) A(g,x)$ \esp for all $x \in X$ and all 
$f,g$ in the acting semigroup $\Gamma$.

\vspace{0.5cm}

\noindent{\bf Theorem A } {\em  Assume that there is a point 
$x_0 \in X$ and a constant $C > 0$ such that for all $f \in \Gamma$, 
$$\max \big\{ \big\| A (f,x_0) \big\|, \big\| A(f, x_0)^{-1} \big\| \big\} \leq C.$$
Then $A$ is cohomologous to a cocycle 
$\tilde{A} \!: \Gamma \times X \to \mathrm{O}(n,\mathbb{R})$, that is, for a 
certain continuous map $B \!: X \to \mathrm{GL}(n,\mathbb{R})$, one has 
\esp $B(f(x))^{-1} A(f,x) B (x) \!=\! \tilde{A}(f,x) \!\in\! O(n,\mathbb{R})$ 
\esp for all $x \!\in\! X$ and all $f \!\in\! \Gamma$.}

\vspace{0.5cm}

This theorem generalizes a classical result of W.H.~Gottschalk and G.A.~Hedlund \cite{GH}, 
which essentially corresponds to the case $n \!=\! 1$. Indeed, Gottschalk and Hedlund 
considered cocycles into the (commutative) group $\mathbb{R}$, which fits in our 
framework by looking at a real number $\lambda$ as the 1-dimensional linear map 
given by multiplication by $e^{\lambda}$; see \S \ref{aplic} for more details.

Thereom A should also be compared with Kalinin's recent 
remarkable extension of Liv$\check{\mathrm{s}}$ic's theorem to matrix 
cocycles \cite{kalinin}. In his setting, the base dynamics is given by that of an Anosov 
diffeomorphisms $T$. Given a H\"older-continuous cocycle $A$ above this dynamics, the 
condition for its cohomological triviality, that is, for the existence of a H\"older-continuous 
$B \!: X \to \mathrm{GL}(n,\mathbb{R})$ such that \esp $A(T,x) = B(Tx) B(x)^{-1}$ \esp 
holds for all $x \in X$, is that the products of $A$ along periodic orbits is trivial:
$$T^n (x) = x \quad \implies \quad \prod_{i=0}^{n-1} A \big( T, T^i (x) \big) = Id.$$  
In view of the method of proof of our Theorem A (see \S \ref{geom}), it is natural 
to ask whether the Kalinin-Liv$\check{\mathrm{s}}$ic theorem admits a version for 
cocycles taking values in the group of isometries of a nonpositively curved space.

\vspace{0.3cm}

\noindent{\bf A criterium of conformality \`a la Sullivan-Tukia.} Let again 
$\Gamma$ be a semigroup acting minimally on a compact metric space $X$, and 
$A \!: X \to \mathrm{GL}(n,\mathbb{R})$ a cocycle above this action. Recall 
that the \textit{quasiconformal distortion} of the linear map $A(f,x)$ is 
defined as 
$$K_A (f,x) := \big\| A (f,x)^{-1} \big\| \cdot \big\| A (f,x) \big\|.$$
Roughly, this measures how distorted is the image under $A(f,x)$ of a 
ball centered at the origin. 

\vspace{0.5cm}

\noindent{\bf Theorem B } If there exists a point $x_0 \in X$ and a constant 
$C > 0$ such that $K_A(f,x_0) \leq C$ holds for all $f \in \Gamma$, then there 
is a continuous invariant conformal structure on the bundle $X\times \mathbb{R}^n$. 
More precisely, the cocycle $A$ is cohomologous to a cocycle taking values in 
the subgroup of conformal linear maps.

\vspace{0.5cm}

This result should be compared with a theorem independently 
due to Sullivan \cite{Su} and Tukia \cite{Tu} (compare also \cite{mark}), 
according to which every uniformly quasiconformal group of diffeomorphisms 
of a 2-manifold is quasiconformally conjugated to a group of conformal maps. Indeed, the first 
step for the proof of this theorem consists in finding an invariant conformal structure; 
the Ahlfors-Bers integrability theorem then allows obtaining the conjugacy. It should be 
pointed out that B. Kalinin and V. Sadovskaya obtained in \cite{KS} an analogous result for 
linear cocycles over an hyperbolic dynamics in the spirit of Liv$\check{\mathrm{s}}$ic's 
theorem.

\vspace{0.5cm}

\noindent{\bf A Bruhat-Tits' lemma in the space of continuous and bounded functions.} 
A very useful lemma due to Bruhat and Tits states that every action by isometries of 
either a proper CAT(0) space or a Hilbert space that has a bounded orbit must have 
a fixed point. Although this still holds for actions on $\mathcal{L}^p$ spaces 
for $1 < p < \infty$, this is no longer true for actions on spaces of continuuous 
functions (see Example \ref{ejem-cont}) and subspaces of $\mathcal{L}^1$ spaces (see 
\cite[Example 2.23]{BFGM}). We next concentrate on the former case in a more 
general situation. 

Let $X$ be a compact metric space, $\mathcal{H}$ a (real) separable Hilbert space, and  
$C(X,\mathcal{H})$ the space of continuous functions on $X$ with values in $\mathcal{H}$. 
In order to discuss affine isometric actions on $C(X,\mathcal{H})$, we need to recall 
a classical result \cite{banach}.

\vspace{0.5cm}

\noindent{\bf Theorem (Banach-Stone)} {\em If $\pi$ is a linear 
surjective isometry of $C(X,\mathbb{R})$, then there exist a unique homeomorphism 
$T\!: X \mapsto X$ and a unique continuous function $sgn \!: X \mapsto \{-1,+1\}$ 
such that for every $\varphi \in C(X,\mathbb{R})$, one has 
$\pi(\varphi)(x) = sgn(x) \varphi ( T^{-1}(x) ).$} 

\vspace{0.5cm}

An almost direct consequence of this theorem is that every action of a 
group $\Gamma$ by linear isometries of $C(X,\mathbb{R})$ comes from an action 
on the basis $X$ together with a cocycle $sgn \!: \Gamma \times X \mapsto \{-1,+1\}$:
$$f\!: \varphi (\cdot) \mapsto sgn \big( f,f^{-1}(\cdot) \big) 
\esp \varphi \big( f^{-1} (\cdot) \big).$$
Here, the cocycle equality is \esp $sgn(fg,x) = sgn (f,g(x)) \esp sgn(g,x)$. 
\esp Moreover, the function $sgn$ must be continuous on the variable $x$.

An analogous statement holds in the space $C(X,\mathcal{H})$ (the corresponding 
version of the Banach-Stone theorem is provided by \cite{jerison}). Thus, every 
action $\pi$ by linear isometries of $C(X,\mathcal{H})$ comes from an 
action (by homeomorphisms) on the basis $X$ together with a cocycle 
$\Psi \!: X \to O(\mathcal{H})$. More precisely, 
\begin{equation} \label{la-rep}
\pi(f) \varphi (x) := \Psi \big( f,f^{-1}(x) \big) \varphi \big( f^{-1}(x) \big), 
\end{equation}
where $\Psi$ satisfies
$$\Psi(fg,x) = \Psi \big( f, g (x) \big) \Psi (g, x).$$

Now let $I \!: \Gamma \to Isom (C(X,\mathcal{H}))$ be an isometric action. By the 
Mazur-Ulam theorem \cite{banach}, $I$ is the composition of a linear isometric action 
$\pi$ and a cocycle $\rho \!: \Gamma \to C(X,\mathcal{H})$, where the cocycle relation is 
$$\rho(fg,x) = \rho(g,x) + \pi(f) \big( \rho(g) \big).$$

\vspace{0.5cm}
 
\noindent{\bf Theorem C } {\em In the context above, assume that the action on the 
basis is minimal. Then the existence of a bounded orbit for the affine isometric action 
\esp $\pi + \rho$ \esp on $C(X,\mathcal{H})$ implies that of a fixed point (function).} 

\vspace{0.5cm}
 
The minimality of the action on $X$ is necessary, as Example \ref{ejem-cont} in \S 
\ref{stat} shows. However, for spaces of bounded measurable functions, there is not 
need to treat any continuity issue, and an analogous (and much simpler~!) version holds 
with no hypothesis on this action. For simplicity, we restrict ourselves to countable 
semigroups (this allows avoiding tedious discussions concerning the measurability of 
certain naturally defined maps). 

\vspace{0.5cm}

\noindent{\bf Theorem D } {\em If an affine isometric action  
$I \!: \Gamma \rightarrow Isom (\mathcal{L}^{\infty}_{\mu}(X,\mathcal{H}))$ 
has a bounded orbit, then it has a fixed point.}

\vspace{0.5cm}

Having our Theorem D as a partial motivation, U. Bader, T. Gelander and N. Monod have 
recently shown an analogous result for $\mathcal{L}^1$ spaces \cite{BGM}. Their clever 
proof is mostly geometric, hence completely different from ours. Quite surprisingly, 
it applies more generally to isometries of preduals of von-Newmann algebras.

Despite the intrinsic interest of Theorems C and D, their possible applications 
in Rigidity Theory are quite limited. Indeed, every countable group acts affinely 
on an $\mathcal{L}^{\infty}$ space and on a space of continuous functions without 
bounded orbits. For instance, one may consider the action of $\Gamma$ on 
$\ell^{\infty}(\Gamma)$ with regular linear part and translation part 
given by \esp $\rho(g)(h) := d(h,g) - d(h,id).$ 


\section{Statement of the Main Theorem and proof of Theorems A, B and C}
\label{stat}

\hspace{0.45cm} As we have already announced, Theorems A, B and C above are almost direct 
consequences of a general principle that is captured by our Main Theorem below. 
Roughly, for every skew action by isometries of a CAT(0) space over a minimal 
dynamics, the existence of a bounded orbit is equivalent to the existence 
of a {\em continuous} invariant section. The proof of Theorem D uses a 
baby form of this principle; see \S \ref{geom}.

Consider a minimal action by continuous maps of a semigroup $\Gamma$ 
on a compact metric space $X$. Let $\mathcal H$ be either a proper 
CAT$(0)$ space or a Hilbert space. We consider a skew action 
by isometries of $\mathcal{H}$:
$$f: (x,v) \mapsto \big( f(x), I(f,x)v \big).$$
Here, for each $f \!\in\! \Gamma$, the map 
$I(f, \cdot) \!: X \rightarrow Isom (\mathcal{H})$ 
is continuous and satisfies the cocycle relation
$$I(fg,x) = I \big( f,g(x) \big) I(g,x).$$ 

\vspace{0.14cm}

\noindent{\bf Main Theorem  } {\em In the setting above, assume that for some $x_0 \in X$ 
and $v_0 \in \mathcal{H}$ there is a bounded subset $B \subset \mathcal{H}$ such that 
$I(f,x_0)v_0$ belongs to $B$ for every $f \in \Gamma$. Then there exists a continuous 
section $x \mapsto (x,\varphi(x)) \in X \times \mathcal{H}$ that is invariant under the 
skew action of $\Gamma$, that is, that satisfies \esp $I(f,x) \varphi(x) = \varphi(f(x))$ 
\esp for all $f \in \Gamma$ and all $x \in X$.}

\vspace{0.5cm}

Notice that the nonpositive curvature hypothesis is necessary, as the 
simple example of an irrational rotation over the torus shows (existence 
of an invariant continuous section is forbidden due to the minimality;  
however, all the orbits are bounded because the underlying product space 
--namely, the torus-- is compact). 

In \S \ref{finito}, we give four independent proofs of the Main Theorem  in the case 
of proper CAT(0) spaces. Letting $\mathcal{H}$ be the hyperbolic plane, this covers 
a case already considered in \cite[Proposition 1]{yoccoz}. (Actually, our fourth 
proof is strongly motivated by that of \cite{yoccoz}.) 

The proof of the Main Theorem for infinite-dimensional Hilbert spaces 
is given in \S \ref{caso-infinito}. This proof is much more subtle than the 
four proofs in \S \ref{finito}. The necessity of a different argument is 
explained by means of a clarifying example in \S \ref{ejem-disc} of a cocycle 
whose linear part is induced by the shift on an orthonormal basis. (These 
cocycles are extensively studied in Appendix B.) Let us 
mention that the argument still applies (with minor modifications that 
we leave to the reader) to the case where the fiber is a (separable) 
uniformly-convex Banach space, thus leading to an analogous theorem in 
this more general situation. The eventual extension to $\mathcal{L}^1$ spaces 
seems to be an interesting problem. Finally, we should point out that, although 
stated for semigroup actions, the Main Theorem extends (with slight modifications 
in the proof) to {\em pseudogroups}, and would also extend to {\em groupoids},  
thus yielding potential applications for foliations.


\vspace{0.1cm}

In what follows, we assume the validity of the Main Theorem, 
and we proceed to give proofs for Theorems A, B and C, together 
with a corollary and an example for the last of these theorems.  

\vspace{0.5cm}

\noindent \textbf{Proof of Theorem A.} The space $Pos(n)$ \esp of positive-definite 
symmetric matrices of order $n \times n$ is a locally-symmetric space of nonpositive 
curvature, hence a proper CAT($0$)-space. The distance between $P \in Pos(n)$ and the 
identity is given by the the sum of the squares of the logarithms of its eigenvalues. 
In particular, there exists $\tilde{C} > 0$ such that 
$\max \{ \| P \|, \| P^{-1} \|\} \leq C$ implies that the distance between 
$P$ and \esp $Id \in Pos(n)$ \esp is smaller than or equal to $\tilde{C}$. 
(See \cite[Chapter XII]{Lang} for more details.)

The linear group $\mathrm{GL}(n,\mathbb{R})$ acts by isometries of $Pos(n)$, 
with $g$ sending $P$ into \esp $g \!\cdot \! P := g P g^{T}.$ The condition 
$\max \{ \|A(f,x_0)\|, \|A(f,x_0)^{-1} \| \} \leq C$ implies that the orbit 
of the point $(x_0,Id)$ under the associated skew action is bounded. By the 
Main Theorem, there exists an invariant continuous section $\varphi \!: X \to Pos(n)$. 

The exponential map at the identity $\exp_{Id}: Sym(n) \to Pos(n)$ is a 
diffeomorphism between the space of symmetric matrices of order $n\times n$ 
and $Pos(n)$. Hence, there is a continuous map $v \!: X \to Sym(n)$ such 
that for each $x\in X,$ 
\[\varphi(x) = \exp_{Id}\left(v(x) \right) \exp_{Id}\left(v(x)\right)^T.\]
We define the continuous map $B \!: X \to \mathrm{GL}(n,\mathbb{R})$ by 
letting $B(x) := \exp_{Id}(v(x))$. (Notice that $B$ takes values in 
$Pos(n)$.) The equation of the invariance of $\varphi$ yields
$$B(f(x))B(f(x))^T
= \varphi(f(x))= A(f,x)\! \cdot \! \varphi(x)
= A(f,x) \varphi(x) A(f,x)^T 
= A(f,x) B(x) B(x)^T A(f,x)^T,$$
hence
$$B(f(x))^{-1}A(f,x) B(x) \big[ B(f(x))^{-1}A(f,x) B(x) \big]^T = Id.$$
Thus, the cocycle $B(f(x))^{-1} A(f,x) B(x)$ takes values in 
$\mathrm{O}(n,\mathbb{R})$, which closes the proof. $\hfill\square$

\vspace{0.5cm}

\noindent \textbf{Proof of Theorem B.} The space \esp $Conf(n)$ \esp of conformal structures 
on $\mathbb{R}^n$ identifies with the space of positive-definite symmetric matrices of order 
$n \times n$ with determinant $1$. This is a Riemannian symmetric subspace of 
$Pos(n)$ with nonpositive curvature, hence a proper CAT(0)-space. The linear group 
$\mathrm{GL}(n,\mathbb{R})$ acts by isometries on $Conf(n)$, with $g$ sending $P$ 
into \esp $g \! \cdot \! P := (\det g^Tg)^{-1/n} g P g^{T}.$ \esp The condition 
$K_A(f,x_0) \leq C$ implies that the orbit of the point $(x_0, Id)$ under the 
associated skew action is bounded. By the Main Theorem, there exists an invariant 
continuous section $\varphi \!: X \to Conf(n)$. As in the proof of Theorem A, 
using the exponential map at the identity, we can find a continuos map 
$B \!: X \to \mathrm{GL}(n,\mathbb R) $ such that for every $x \in X$,
$$\varphi (x) = B(x)B(x)^T.$$
Denote $\lambda(x)=(\det A(f,x)^TA(f,x))^{-1/2n}$. 
The invariance of $\varphi$ yields
\begin{eqnarray*}
B(f(x))B(f(x))^T \hspace{0.2cm} = \hspace{0.2cm} \varphi (f(x)) \hspace{0.2cm} = 
\hspace{0.2cm} A(f,x) \cdot \varphi(x) &=& \lambda(x)^2 A(f,x) \varphi (x) A(f,x)^T\\ 
&=& \lambda(x)A(f,x) B(x) \big[ \lambda(x)A(f,x) B(x) \big]^T,
\end{eqnarray*}
hence
$$\lambda(x) B(f(x))^{-1} A(f,x) B(x) \big[ \lambda(x) B(f(x))^{-1} A(f,x) B(x) \big]^T = Id.$$
We thus conclude that the cocycle $\lambda(x) B(f(x))^{-1} A(f,x) B(x)$ takes 
values in $\mathrm{O}(n,\mathbb{R})$, and therefore $B(f(x))^{-1} A(f,x) B(x)$ 
belongs to the conformal linear group of $\mathbb{R}^n$. $\hfill\square$

\begin{obs} 
{Notice that if it is possible to solve the classical cohomological equation for $\lambda$, 
then this allows conjugating $A$ into a cocycle taking values in $\mathrm O(n, \mathbb R)$.}
\end{obs}

\vspace{0.3cm}

\noindent \textbf{Proof of Theorem C.} Writing  \esp $\rho(f,x) := \rho (f) (f^{-1}(x))$, 
\esp so that the isometric action may be written as 
$$I(f) \varphi(x) = \Psi \big( f,f^{-1}(x) \big) \varphi \big( f^{-1}(x) \big) 
+ \rho \big( f,f^{-1}(x) \big),$$
we have the cocycle relations 
$$\Psi(fg,x) = \Psi \big( f,g(x) \big) \Psi(g,x), \qquad 
\rho(fg,x) = \rho \big( f,g(x) \big) + \Psi \big( f,g(x) \big) \rho(g,x).$$
It is then easy to check that  
$$f \!: (x,v) \mapsto \Psi (f,x) v + \rho(f,x)$$ 
defines a skew action on $X \times \mathcal{H}$ by isometries on the fibers. Since $I$ 
is assumed to have a bounded orbit, all its orbits must be bounded. In particular, the 
orbit of the identically zero function is bounded, that is, there exists a constant 
$C$ such that $\| \rho(f,f^{-1}(x)) \| \leq C$ holds for all 
$f \in \Gamma$ and all $x \in X$. This means 
that the orbit of the zero vector of $\mathcal{H}$ under the associated skew action 
on $X \times \mathcal{H}$ is bounded. By Theorem A, there exists a continuous 
function $\varphi_0 \!: X \to \mathcal{H}$ satisfying, for all $f \in \Gamma$ 
and all $x \in X$,
$$\varphi_0 \big( f(x) \big) = \Psi(f,x) \varphi_0 (x) + \rho(f,x).$$
Changing $x$ by $f^{-1}(x)$, this equality becomes
$$\varphi_0 (x) = \Psi \big( f,f^{-1}(x) \big) \varphi_0 \big( f^{-1}(x) \big) 
+ \rho(f,f^{-1}(x)) = I(f) \varphi_0 (x),$$   
thus showing that $\varphi_0 \in C(X,\mathcal{H})$ is a fixed point of $I$. $\hfill\square$

\vspace{0.55cm}

The next corollary to Theorem C was kindly suggested to the second-named 
author by P. Py, and should be compared with the results of Appendix A.

\vspace{0.35cm}

\noindent{\bf Corollary} {\em Consider a linear representation $\pi$ on 
$C(X,\mathcal{H})$ of the form} (\ref{la-rep}). {\em If the $\Gamma$-action 
on $X$ is minimal, then for every quasi-invariant probability measure $\mu$ 
on $X$, the natural map from $H^1 \big( \pi,C(X,\mathcal{H}) \big)$ into 
$H^1 \big( \pi,\mathcal{L}^{\infty}_{\mu}(X,\mathcal{H}) \big)$ is injective.}

\vspace{0.25cm}

\noindent{\bf Proof.} Let $\rho \!: \Gamma \to C(X,\mathcal{H})$ be a cocycle 
that is cohomologically trivial in $\mathcal{L}^{\infty}_{\mu}(X,\mathcal{H})$. 
Due to Theorem C, we need to show that $\rho(g)$ is bounded as a function 
in $C(X,\mathcal{H})$ independently of $g$. To do this, we may assume 
that $\Gamma$ is countable. Indeed, if $\rho$ is not bounded, then there exists 
a sequence $g_n \in \Gamma$ such that $\|\rho(g_n)\|_{C(X,\mathcal{H})} \geq n$, 
for each $n \in \mathbb{N}$. Thus, the cocycle $\rho$ is unbounded when restricted 
to the countably generated subgroup $\langle g_1,g_2,\ldots \rangle$.

Now, since $\rho$ is cohomologically trivial in 
$\mathcal{L}^{\infty}_{\mu}(X,\mathcal{H})$, it may be written in the form 
\begin{equation}\label{casi}
\rho \big( g,g^{-1}(x) \big) = \rho(g)(x) 
= \Psi \big( g,g^{-1}(x) \big) \varphi \big( g^{-1}(x) \big) - \varphi(x)
\end{equation}
for a certain function $\varphi \in \mathcal{L}^{\infty}_{\mu}(X,\mathcal{H})$, 
where the second equality above holds $\mu$-a.e. Let $X_0$ be the set of points 
$x \in X$ for which equality (\ref{casi}) does not hold for some $g \in \Gamma$. 
Since $\Gamma$ is assumed to be countable, $X_0$ has zero $\mu$-measure. Let 
$C$ be the essential supremum of $\| \varphi \|$. Then the $\mu$-measure of 
$X^* := \{x \!: \| \varphi(x) \| > C \}$ is zero, as well as that of 
$X_1 := \bigcup_{g \in \Gamma} g^{-1} (X^*)$. Let $x_0$ be a point 
in the full $\mu$-measure set $X \setminus (X_0 \cup X_1)$. Then 
equality (\ref{casi}) holds at $x_0$ for all $g \in \Gamma$. 
Moreover, \esp $\| \varphi \big( g(x_0) \big) \| \leq C$ 
also holds for all $g \in \Gamma$. 
This allows us to conclude that, for all $g \in \Gamma$, 
\begin{equation}\label{first-estimate}
\big\| \rho \big( g,g^{-1}(x_0) \big) \big\| \leq 2C. 
\end{equation}
We claim that for all $h \in \Gamma$, we have $\| \rho(h) \| \leq 4C$. 
Indeed, the cocycle identity yields
$$\rho \big( gh, (gh)^{-1}(x_0) \big) = 
\Psi \big( g,g^{-1}(x_0) \big) \esp \rho \big( h, (gh)^{-1} (x_0 ) \big) + 
\rho \big( g,g^{-1}(x_0) \big).$$
Thus, by (\ref{first-estimate}), 
$$\big\| \rho \big( h,h^{-1}g^{-1}(x_0) \big) \big\|
\leq \big\| \rho \big( gh,(gh)^{-1} (x_0) \big) \big\| 
+ \big\| \rho \big( g,g^{-1}(x_0) \big) \big\| \leq 4C.$$
Fix $x \in X$. Taking a sequence $(g_n)$ in $\Gamma$ such that $g_n^{-1} (x_0) \to x$ 
as $n \to \infty$, we obtain  
$$\big\| \rho(h)(x) \big\| = \big\| \rho \big( h,h^{-1}(x) \big) \big\| = 
\lim_{n \to \infty} \big\| \rho \big( h, h^{-1} g_n^{-1} (x_0) \big) \big\| 
\leq 4C,$$
which shows our claim and hence the Corollary.  $\hfill\square$

\vspace{0.32cm}

We close this section with an example showing that the hypothesis of minimality for 
the action on $X$ above is necessary. (A more interesting example in that the action 
on the basis is topologically transitive can be derived from \cite[Exercise 2.9.2]{KH}.)

\vspace{0.1cm}

\begin{ejem} \label{ejem-cont}
{\em Consider a parabolic element $T \in \mathrm{PSL}(2,\mathbb{R})$ 
acting on $X := \mathrm{S}^1$. Denoting by $x_0$ the unique fixed point of $T$, 
we let $\psi \!: \mathrm{S}^1 \mapsto \mathbb{R}$ be a function having a single 
discontinuity at $x_0$, so that $T(x_0)$ equals \esp $\lim_{x \to x_0^+} T(x)$ 
\esp and is different from \esp $\lim_{x \to x_0^{-}} T(x)$. \esp Then the function 
$x \mapsto \psi - \psi \circ T$ is continuous (it vanishes at $x_0$). Therefore, 
we may consider the affine isometric action of $\Gamma \sim \mathbb{Z}$ on 
$C(\mathrm{S}^1,\mathbb{R})$ generated by 
$$I(T) \varphi (x) := \varphi \big( T^{-1}(x) \big) + \psi(x) - \psi \big( T^{-1}(x) \big).$$
Since, for every $n \in \mathbb{Z}$, 
$$I(T^n) \varphi (x) = \varphi \big( T^{-n}(x) \big) + \psi(x) - \psi \big( T^{-n}(x) \big),$$
the orbit of any $\varphi \in C(\mathrm{S}^1,\mathbb{R})$ is bounded in norm 
by \esp $\| \varphi \|_{C(X,\mathbb{R})} + 2 \|\psi\|_{\mathcal{L}^{\infty}}$. 
\esp We claim that, however, there is no fixed point in $C(\mathrm{S}^1,\mathbb{R})$ 
for this action, so that the cocycle  \esp $\psi - \psi \circ T$ \esp \esp is  
trivial in $H^1 \big( \pi,\mathcal{L}^{\infty}_{\mu}(X,\mathcal{H}) \big)$ 
but nontrivial in $H^1 \big( \pi,C(X,\mathcal{H}) \big)$. Indeed, 
the equality $I(T^{-1}) \varphi = \varphi$ yields, for 
every $x \in \mathrm{S}^1$ and all $n \in \mathbb{N}$, 
$$\varphi - \psi = (\varphi - \psi) \circ T = \ldots = 
(\varphi - \psi) \circ T^n.$$
Since the (forward) $T$-orbit of any $x \in \mathrm{S}^1$ 
converges to $x_0$, say by the right, this implies that the value of 
$\varphi \!-\! \psi$ is constant and equals \esp $\varphi(x_0) \!-\! 
\lim_{x \to x_0^+} \psi(x)$. \esp Clearly, this implies that 
$\varphi$ cannot be continuous.} 
\end{ejem}


\section{Further applications: cohomological equations}
\label{aplic}

\hspace{0.45cm} Several problems in dynamical systems reduce to solving a linear functional 
(or {\em cohomological}) equation. For example, the (linearized version of the) conjugacy 
problem for circle diffeomorphism (see \cite{herman}), the study of interval exchange maps 
(see \cite{MMY}), the existence of eigenvalues of the Koopman operator associated with 
a dynamical system (see \cite{KARO}), time changes for flows (see \cite{KH}), etc. One 
of the most basic results about the existence of continuous solutions for these 
equations is the classical Gottschalk-Hedlund theorem that we next recall (see 
\cite[Chapter 14]{GH} for more details). Notice that the converse of this result 
is also true but much more elementary. 
  
\vspace{0.4cm}

\noindent{\bf Theorem (Gottschalk-Hedlund)} {\em  Let $X$ be a compact metric space, 
$T:X\to X$ a minimal continuous map and $\rho:X\to \R$ a continuous function. If 
there exists a point $x_0\in X$ such that
\begin{equation}\label{bd}
\sup_{n\in \N}\Big |\sum_{j=0}^{n-1}\rho \big( T^j (x_0) \big) \Big|<\infty,
\end{equation}
then the cohomological equation
\begin{equation}\label{EWP-P}
\varphi\circ T - \varphi =\rho
\end{equation}
has a continuous solution $\varphi: X \to \R$.}

\vspace{0.4cm}

The origin of the Gottschalk-Hedlund theorem was the study of a special 2-dimensional system, 
nowadays known as {\em cylindrical cascade}. Let $X, T$ and $\rho$ be as before. The 
cylindrical cascade associated to this data is the map
\begin{eqnarray*}
F: X\times \R&\to& X\times \R\\
(x, t)&\mapsto& \big( T(x), t+\rho(x) \big).
\end{eqnarray*}
Gottschalk and Hedlund observed that $F$ is topologically conjugated to 
the map $(x, t) \mapsto (T(x), t)$ if and only if the cohomological 
equation (\ref{EWP-P}) has a continuous solution. 

The map $F$ above can be though of as the skew action induced by a minimal $\mathbb{N}$-action 
on $X$ and a cocycle of isometries (translations) of $\R$. Moreover, the hypothesis (\ref{bd}) 
corresponds to that the orbit of the point $(x_0,0)$ under this skew action is bounded. This 
fits into both the framework and the hypothesis of our Main Theorem for the case of a cocycles 
$I$ into the group of isometries of a Hilbert space 
$\mathcal{H}$. Indeed, writing $I = \Psi + \rho$, with 
$\Psi$ being the linear part of $I$ and $\rho$ being the translation part, the cocycle 
relations become 
\begin{equation}\label{cociclos-GH}
\Psi(fg,x) = \Psi \big( f,g(x) \big) \Psi(g,x), \qquad 
\rho(fg,x) = \rho \big( f,g(x) \big) + \Psi \big( f,g(x) \big) \rho(g,x).
\end{equation}
Whenever this is satisfied, we have an associated skew action on $X \times \mathcal{H}$:
$$f: (x,v) \mapsto \big( f(x), I(f,x)v \big).$$ 
The Main Theorem asserts that the existence of a bounded orbit for this skew 
action implies the existence of a continuous invariant section $\varphi$.
Since this means that \hspace{0.01cm}
$I(f,x) \varphi(x) = \varphi(f(x)),$ \hspace{0.01cm}
we have that $\varphi$ satisfies the twisted cohomological equation
\begin{equation}\label{chupalo-vos}
\varphi \big( f(x) \big) - \Psi(f,x) \varphi(x) = \rho(f,x).
\end{equation}
Moreover, conjugation by the homeomorphism 
$S\!: (x,v) \mapsto \big( x,v-\varphi(x) \big)$ yields, for each $f \in \Gamma$,
\begin{eqnarray*}
S f S^{-1} (x,v) 
&=& S f \big( x, v + \varphi(x) \big)\\ 
&=& S \big( f(x), I(f,x)(v + \varphi(x)) \big)\\ 
&=& S \big( f(x), \Psi(f,x)v + I(f,x) \varphi(x) \big)\\ 
&=& \big( f(x), \Psi(f,x)v + I(f,x) \varphi(x) - \varphi(f(x)) \big)\\ 
&=& \big( f(x), \Psi(f,x)v \big).
\end{eqnarray*}
In other words, conjugation by $S$ reduces the cocycle $I$ to its linear part $\Psi$.

\begin{ejem} \label{GH}
{\em If $\mathcal{H} = \mathbb{R}$ and $\Psi(f,x) = Id$ for every $(f,x)$, then the Main Theorem  
is the version for semigroups of the ``equivariant Gottschalk-Hedlund lemma'' of \cite{yo} 
(also contained in \cite[Section 3.6.2]{libro}), which -\!~as the second-named author 
discovered while writing this article~\!- was originally obtained by J.~Moulin Ollagnier 
and D.~Pinchon in \cite{el-primero} (compare \cite{KL}). Notice that for 
$\Gamma \sim \mathbb{N}$, this corresponds to the classical Gottschalk-Hedlund 
theorem. Nevertheless, even in this particular case, the proof we will provide 
for the Main Theorem  differs from the classical ones in a key geometric argument. 
For $\Gamma \sim \mathbb{R}^+$, this is an equivalent form of the main 
result of \cite{mc}.}
\end{ejem}

\begin{ejem} \label{plus-tard}
{\em Again in dimension 1, let $\Gamma \sim \mathbb{N}$ act on $X$ 
by powers of a continuous, minimal map $T$. Letting $\Psi(n,x) := (-Id)^n$, 
the Main Theorem  yields the following statement: if, for a continuous function 
$\rho \!: X \rightarrow \mathbb{R}$, the values of the alternating sums 
$$\sum_{k=0}^{n-1} (-1)^{k} \rho \big( T^k(x_0) \big),$$
are uniformly bounded (independently of $n$) for some $x_0 \in X$, 
then the cohomological equation
$$\varphi \big( T(x) \big) + \varphi(x) = \rho(x)$$
has a continuous solution $\varphi$. The interest on this equation comes 
from the problem of extracting a square root of the associated cylindrical 
cascade. More precisely, if $T^{1/2}$ is a square root of $T$, then 
$(x,v) \mapsto ( T^{1/2}(x), v + \varphi(x) )$ is a square root 
of $(x,v) \mapsto ( T(x), v + \rho(x) )$ if and only of $\varphi$ 
satisfies the cohomological equation}
$$\varphi \big( T^{1/2}(x) \big) + \varphi(x) = \rho(x).$$ 
\end{ejem}

\begin{ejem} \label{ejem-melnikov} 
{\em Let $\mathcal{H} = \R^2 \sim \mathbb{C}$, 
and consider an action of $\Gamma \sim \mathbb{N}$ on $X$ by 
powers of a continuous, minimal map $T$. Assume that $\Psi(n,x) = \Psi(n)$ 
does not depend on $x$ and preserves orientation. Then it coincides with 
the rotation of angle $n \beta$, where $e^{i\beta} = \Psi(1,x)$ for any 
$x$. Given $z \in \mathbb{C}$, the cocycle relation yields
$$I(n,x) z = \Psi(n) z + \rho(n,x) 
= e^{in\beta} z + \sum_{k=0}^{n-1} e^{i(n-k-1)\beta} \rho \big( 1,T^{k}(x) \big).$$
In this case, the boundedness hypothesis means that for 
$\rho(x) := \rho(1,x)$ and some $x_0 \in X$, the norm of
\begin{equation}\label{gatito1}
\sum_{k=0}^{n-1} e^{-i k \beta} \rho \big( T^{k}(x_0) \big)
\end{equation}
is uniformly bounded (independently of $n$). Moreover,
equation (\ref{chupalo-vos}) becomes
\begin{equation}\label{gatito2}
\varphi \big( T(x) \big) - e^{i \beta} \varphi(x) = \rho(x).
\end{equation}
As a consequence of the Main Theorem, if the sums (\ref{gatito1}) 
are uniformly bounded, then the cohomological equation (\ref{gatito2}) 
has a continuous solution $\varphi$. (An alternative, less geometric   
proof of this fact can be derived from the results of \cite{ponce}.) 
We should point out that this equation corresponds to the linearized version 
of that encoding the stability of a closed orbit under perturbation in a toy model 
of a planetary system; see \cite{poschel}. We refer to \cite{CNP-london} for an 
accurate study of the dynamics of the associated map 
$(x,z) \mapsto (T(x),e^{i\beta}z + \rho(x))$.}
\end{ejem}

\begin{ejem} \label{ejem-ciclon}
{\em Given an integer $q \geq 1$, an irrational angle $\alpha$, an arbitrary 
angle $\beta$, and a continuous function $\rho \!: \clo \to \mathbb{C}$, we consider 
the skew map $(\theta, z) \mapsto (\theta + \alpha, e^{i\beta}z + \rho(\theta) )$ from 
$\clo \times \mathbb{C}$ into itself. (This corresponds to a particular case of Example 
\ref{ejem-melnikov}.) One easily checks that finding a $q^{th}$ root of the form 
$(\theta,z) \mapsto (\theta + \alpha/q, e^{i\beta / q} z + \varphi(\theta) )$ 
for this map is equivalent to solving the {\em cyclotonic equation}} 
\begin{equation}\label{ciclon-eq}
\sum_{k=0}^{q-1} e^{\frac{ik \beta}{q}} \varphi \Big( \theta + \frac{(q-k-1)\alpha}{q} \Big) 
= \rho (\theta).
\end{equation}
{\em Notice that, for $\beta := 0$ and $q := 2$, we retrieve an equation similar to 
that of Example \ref{plus-tard}.

Despite the strange form of equation (\ref{ciclon-eq}), we claim that if $\alpha$ 
and $\beta$ are independent over the rationals, then it is equivalent to the 
twisted cohomological equation}
\begin{equation}\label{reduction}
\varphi(\theta + \alpha) - e^{i \beta} \varphi(\theta) = 
\rho \Big( \theta + \frac{\alpha}{q} \Big) - 
e^{\frac{i\beta}{q}} \rho \Big( \theta + \frac{\alpha}{q} \Big).
\end{equation}

{\em Indeed, if $\varphi$ solves (\ref{ciclon-eq}), then changing $\theta$ by 
$\theta - \alpha/q$ and multiplying both sides by $e^{\frac{i\beta}{q}}$, we obtain 
$$\sum_{k=0}^{q-1} e^{\frac{i(k+1)\beta}{q}} 
\varphi \Big( \theta + \frac{(q-1-(k+1)) \alpha}{q} \Big) 
= e^{\frac{i \beta}{q}} \rho \Big( \theta - \frac{\alpha}{q} \Big),$$
that is} 
\begin{equation}\label{despl}
\sum_{k=1}^{q} e^{\frac{i k \beta}{q}} \varphi \Big( \theta + \frac{(q-k-1) \alpha}{q} \Big) 
= e^{\frac{i\beta}{q}} \rho \Big( \theta - \frac{\alpha}{q} \Big).
\end{equation}
{\em Substracting (\ref{despl}) from (\ref{ciclon-eq}) yields} 
$$\varphi \Big( \theta + \frac{(q-1) \alpha}{q} \Big) - 
e^{\frac{iq\beta}{q}} \varphi \Big( \theta - \frac{\alpha}{q} \Big) = 
\rho(\theta) - e^{\frac{i\beta}{q}} \rho \Big( \theta - \frac{\alpha}{q} \Big).$$
{\em Finally, changing $\theta$ by \esp $\theta + \alpha/q$ \esp yields (\ref{reduction}).

Conversely, assume that $\varphi$ solves (\ref{reduction}). Then} 
\begin{eqnarray*}
\sum_{k=0}^{q-1} e^{\frac{ik\beta}{q}} \varphi \Big( \theta + \frac{(q-k-1) \alpha}{q} \Big) 
\!\!\!&=&\!\!\! \sum_{k=0}^{q-1} e^{\frac{ik\beta}{q}} 
     \left[ e^{i\beta} \varphi \Big( \theta - \frac{(k+1)\alpha}{q} \Big)  
    \! + \rho\Big( \theta - \frac{k\alpha}{q}\Big) \!- 
      e^{\frac{i\beta}{q}} \rho \Big(\theta - \frac{(k+1)\alpha}{q}\Big) \! \right]\\
&=&\!\!\! \rho (\theta) - e^{i\beta} \rho(\theta - \alpha) + 
    e^{i\beta} \sum_{k=0}^{q-1} e^{\frac{ik\beta}{q}} 
         \varphi \Big( \theta - \frac{(k+1)\alpha}{q} \Big).
\end{eqnarray*}
{\em Letting} 
$$\psi(\theta) := 
\sum_{k=0}^{q-1} e^{\frac{ik\beta}{q}} \varphi \Big( \theta + \frac{(q-k-1) \alpha}{q} \Big),$$
{\em this equality may be rewritten as}
$$\psi(\theta) = 
e^{i\beta} \psi(\theta - \alpha) + \rho(\theta) - e^{i\beta} \rho(\theta - \alpha),$$
{\em that is,}
$$\psi(\theta) - \rho(\theta) 
= e^{i\beta} \big[ \psi(\theta - \alpha) - \rho(\theta - \alpha) \big].$$
{\em The function $\theta \mapsto \| \psi (\theta) - \rho(\theta) \|$ is thus invariant under 
the rotation of angle $-\alpha$, hence constant. If the value of this constant is nonzero, then 
it is well known that (modulo 1) $\beta$ must be a rational multiple of $\alpha$ (see, for 
instance, \cite[Theorem 3.5]{walters}), which is contrary to our hypothesis.}
\end{ejem}


\section{Proof of the Main Theorem}

\subsection{A general strategy and proof of Theorem D}
\label{geom}

\hspace{0.45cm} An important case covered by the Main Theorem corresponds to that where $X$ 
is a single point. In this case, our result reduces to the version for semigroups of the 
Bruhat-Tits lemma. To better discuss this link, we recall the general framework (see 
\cite[Proposition 5.10]{ballmann} for more details). Let $\mathcal{H}$ be either a proper 
CAT(0) space or a (real and separable) Hilbert space. Given a bounded subset of 
$\mathcal{H}$, for each $v \in \mathcal{H}$ we let 
$$r_B (v) := \inf \big\{r > 0 \!: B \subset 
\mathrm{Ball}(v,r) \big\} = \sup_{w \in B} d(v , w ).$$ 
The {\em radius} of $B$ is defined as \hspace{0.01cm} 
$r_B := \inf \big\{ r_B (v) \!: v \in \mathcal{H} \big\}.$ 
\hspace{0.01cm} The following facts hold:

\vspace{0.1cm}

\noindent -- The infimum of $r_B (\cdot)$ is attained. Indeed, in case of a proper 
CAT(0) space, this follows from the compactness of the closed (bounded) balls. In 
case of a Hilbert space, this follows from the relative compactness of bounded 
subsets of $\mathcal{H}$ when endowed with the weak topology, and the fact that 
the distance function is lower-semicontinuous.

\vspace{0.1cm}

\noindent -- Actually, it is attained at a unique point. Inded, this follows from 
the ``convexity properties'' of the distance function on $\mathcal{H}$, that is, 
the CAT(0) property.

\vspace{0.1cm}

\noindent The unique point realizing the infimum is called the geometric (or Chebyshev) 
{\em center} of $B$. This point $w := ctr(B)$ is thus characterized as being the unique 
one satisfying $B \subset \overline{\B (w,r_B)}$. 

By construction, if $I\!: \mathcal{H} \to \mathcal{H}$ is an isometry, 
then $r_B = r_{I(B)}$ and $I(ctr(B)) = ctr (I(B))$. Moreover, we have 
the following fact (a proof is given further on).

\vspace{0.05cm}

\begin{prop} \label{se-usa} The map $B \mapsto ctr(B)$ is continuous with respect
to the Hausdorff topology on bounded subsets of $\mathcal{H}$.
\end{prop}

\vspace{0.05cm}

Let us again recall the statement of the Bruhat-Tits center lemma \cite{BT}.

\vspace{0.45cm}

\noindent{\bf Lemma (Bruhat-Tits)} {\em Let $\Gamma$ be a group acting by isometries 
of $\mathcal{H}$. If the action has a bounded orbit, then there is a point in 
$\mathcal{H}$ that is fixed by every element of $\Gamma$. As a consequence, 
the action of $\Gamma$ is conjugate to an action by linear isometries.}

\vspace{0.45cm}

Indeed, the center of the bounded orbit must remain fixed. 
If we conjugate by the translation sending this fixed
point to the origin, then the action of every
element of $\Gamma$ becomes an isometry fixing
the origin, that is, a linear isometry.

It is worth mentioning that Bruhat-Tits' lemma still holds for {\em semigroup} actions, but the 
proof needs an extra argument. Indeed, if $B$ is a bounded forward-invariant set (as for example
a bounded orbit of the semigroup), it is not completely obvious that its center is invariant 
by every $f \!\in\! \Gamma$. To see that this is the case, notice that, letting $r:= r_B$, 
from $B \subset \overline{\mathrm{Ball} (ctr(B),r)}$ we obtain 
$f(B) \subset \overline{\mathrm{Ball} (f(ctr(B)),r)}$. 
Now, as $f (B) \subset B$, we also have $f(B) \subset \overline{\B (ctr(B),r)}$. Since 
$r = r_{f(B)}$, this necessarily implies that $f(ctr(B)) = ctr(B)$, as desired.

\vspace{0.5cm}

\noindent{\bf The main idea.} The construction above provides us with a basic strategy 
of proof for the Main Theorem. Indeed, according to the hypothesis, the $\Gamma$-orbit 
of certain point $(x_0,v_0) \in X \times \mathcal{H}$ 
remains in a bounded subset of $X \times \mathcal{H}$. By continuity, 
its closure $M := \overline{orb(x_0,v)}$ is a compact, forward-invariant set.
Notice that since the $\Gamma$-action on $X$ is assumed to be minimal, the projection of
$M$ on $X$ is the whole space. As a consequence, the $\Gamma$-orbit of {\em any} point 
$(x,v)$ remains in a bounded subset of $X \times \mathcal{H}$ (which depends on $(x,v)$). 
Indeed, if $v^* \in \mathcal{H}$ is such that $(x,v^*)$ belongs to $M$, then the 
$\Gamma$-orbit of $(x,v^*)$ is contained in $M$. Since for each $v \in \mathcal{H}$ 
and all $f \in \Gamma$,
$$d \big( I(f,x) v , I(f,x) v^* \big) = d(v , v^* ),$$
this implies that the $\Gamma$-orbit of $(x,v)$ is also bounded.

For each $x \in X$, let $M_x := \{v \in \mathcal{H} \!: (x,v) \in M \}$. 
Notice that $I(f,x) M_x = M_{f(x)}$ holds for all $f \in \Gamma$ and all 
$x \in X$. The (nonempty) set $M_x$ is bounded, hence we may consider 
its center $\varphi(x) := ctr (M_x)$. Since the center map commutes with 
isometries, the curve $x \mapsto (x,\varphi(x))$ is invariant under the skew 
action. However, it is not evident at all that the thus-obtained map $\varphi$ 
is continuous ({\em a priori}, it is just {\em measurable}). Indeed, we will need to 
elaborate a little bit to show that this is always the case for proper spaces. For  
infinite-dimensional Hilbert space fibers, this may fail to happen, hence we will 
need to slightly modify our approach. The proof for this case is strongly motivated 
by the main argument of Namioka-Asplund's proof of the Ryll-Nardewski fixed point 
theorem \cite{NA}. Let us point out that a slight modification allows applying this 
argument also for CAT(0) proper spaces. More importantly, it easily applies to 
cocycles of noncontracting maps of $\mathcal{H}$, thus extending our Main 
Theorem to this framework.

\vspace{0.35cm}

As a first illustration of the preceding idea, we next give a 

\vspace{0.43cm}

\noindent{\bf Proof of Theorem D.} The proof is similar to that of Theorem C 
though much simpler since we do not need to take care of continuity issues. 
Let $\rho \!: \Gamma \to \mathcal{L}^{\infty}_{\mu}(X,\mathcal{H})$ be 
the translation part associated to the representation $I$, so that 
\begin{equation}\label{action-inf}
I(f) \varphi (x) 
= \Psi \big( f,f^{-1}(x) \big) \varphi \big( f^{-1}(x) \big) + \rho \big( f,f^{-1}(x) \big).
\end{equation}
Assume that the $I$-orbit $orb(\varphi_0)$ of $\varphi_0$ is bounded so that 
the norm of each point therein is less than or equal to 
a constant $C$. For each $x \in X$, we let \esp 
$N_x := \{\varphi(x) \!: \varphi \in orb(\varphi_0) \}.$ \esp 
Then for $\mu$-almost-every $x \in X$, this set $N_x$ is bounded in norm 
by $C$. We may thus consider the function $\varphi \!: X \to \mathcal{H}$ 
defined by $\varphi(x) := ctr (N_x)$. One can check that this is a 
measurable function. (This is an easy exercise if $\mathcal{H}$ has finite 
dimension, but a little bit harder in the infinite-dimensional case.) 
Moreover, it clearly belongs to $\mathcal{L}^{\infty}_{\mu}(X,\mathcal{H})$. 
We claim that $\varphi$ is a fixed point of the isometric action. Indeed, 
due to (\ref{action-inf}), for $\mu$-almost-every $x \in X$ and 
every $g \in \Gamma$, we have 
\begin{eqnarray*}
N_{g^{-1}(x)} 
&=& \big\{\varphi(g^{-1}(x)) \!: \varphi \in orb(\varphi_0) \big\}\\ 
&=& \big\{ \Psi \big( g,g^{-1}(x) \big)^{-1} \big[ I(g) \varphi (x) - \rho \big(g,g^{-1}(x)\big) \big] 
\!: \varphi \in orb(\varphi_0) \big\}\\ 
&=& \Psi \big( g,g^{-1}(x) \big)^{-1} \big\{ I(g) \varphi (x)  \!: \varphi \in orb(\varphi_0) \big\} 
    - \Psi \big( g,g^{-1}(x) \big)^{-1} \rho \big( g, g^{-1}(x) \big)\\ 
&=& \Psi \big( g,g^{-1}(x) \big)^{-1} (N_x) 
    - \Psi \big( g,g^{-1}(x)\big)^{-1} \rho \big( g,g^{-1}(x) \big), 
\end{eqnarray*}
hence 
$$\Psi \big( g, g^{-1}(x) \big) (N_{g^{-1}(x)}) + \rho \big( g, g^{-1}(x) \big) =  N_x.$$
Taking the center at both sides we obtain, for $\mu$-almost-every $x \in X$,
$$\Psi \big( g, g^{-1}(x) \big) \varphi \big( g^{-1}(x) \big) 
+ \rho \big( g, g^{-1} (x) \big) = \varphi(x),$$
which is equivalent to $I(g)\varphi = \varphi$. $\hfill\square$


\subsection{The finite-dimensional case}
\label{finito}

\subsubsection{First proof}

\hspace{0.45cm} In the context of the Main Theorem, assume that 
$\mathcal{H}$ is a proper CAT(0) space. Let $M$ be any nonempty, 
compact invariant set for the skew action on $X \times \mathcal{H}$.

\vspace{0.05cm}

\begin{lema} \label{unila}
Given $x \in X$, let $(f_k)$ be a sequence of elements in $\Gamma$
such that $f_k(x) \rightarrow x$. Then $I(f_k,x) (M_x)$ converges
(in the Hausdorff topology) to $M_x$.
\end{lema}

\noindent{\bf Proof.} If not, then one of the 
following two possibilities should arise.

\vspace{0.2cm}

\noindent {\bf 1.-} There is a sequence of points 
$v_k \in I(f_k,x)(M_x)$ converging to a certain $v^* \notin M_x$.

\vspace{0.1cm}

This case is impossible. Otherwise, the sequence of points $(f_k(x),v_k) \in M$
would converge to the point $(x,v^*) \notin M$, thus contradicting the fact that
$M$ is closed.

\vspace{0.2cm}

\noindent {\bf 2.-} There is a point $v^* \in M_x$
having a neighborhood $V \subset \mathbb{R}^{\ell}$ such that,
for large-enough $k$, no point $w \in I(f_k,x) (M_x)$ belongs to $V$.

\vspace{0.1cm}

This case is impossible as well, but the argument is more subtle. First, notice that 
since $I(f_k,x)(M_x)$ is uniformly bounded on $k$, the isometries $I(f_k,x)$ remain
inside a compact subset of $Isom (\mathcal{H})$. Passing to a subsequence if
necessary, we may assume that they converge to some $I \in Isom (\mathcal{H})$.
By our assumption, the vector $v^*$ belongs to $M_x \setminus I(M_x)$. Moreover,
we must have $I(M_x) \subset M_x$, because $M$ is invariant and closed. Therefore,
$I(M_x) \subsetneq M_x$. Since $M_x$ is closed, this is impossible, because of
the following

\vspace{0.2cm}

\noindent{\underbar{Independent Claim.}} If $C$ is a (nonempty) compact subset of 
$\mathcal{H}$ and $J$ is an isometry such that $J(C) \subset C$, then $J(C) = C$.

\vspace{0.1cm}

To show this, first notice that $J^n (C) \subset C$ for all $n \in \mathbb{N}$, which
forces the subgroup $G = \overline{\langle J \rangle}$ of $Isom(\mathcal{H})$ to be
compact. We claim that $Id$ is an accumulation point of $G$.
Indeed, if $J_0$ is any accumulation point
of $G$, then for every $\varepsilon > 0$, there exists $n_1 < n_2$ such that
$dist (J^{n_{1}},J_0) < \varepsilon / 2$ and
$dist (J^{n_{2}},J_0) < \varepsilon / 2$, where $dist$ is a bi-invariant 
metric on $G$. This yields $dist (J^{n_{1} - n_{2}}, Id) < \varepsilon$. Since
this holds for every $\varepsilon > 0$, there exists an increasing sequence of
integers $m_k$ such that $J^{m_k}$ converges to the identity, as asserted.

Assume now that some vector $v$ belongs to $C \setminus J(C)$. Since $J(C)$ is 
closed, there exists a neighborhood $V$ of $v$ such that $J(C) \cap V = \emptyset$. 
Then we have $J^m (C) \cap V = \emptyset$ for all $m > 0$. However, since 
$J^{m_k} \rightarrow Id$, the sets $J^{m_k}(C)$ converge to $C$ in the 
Hausdorff topology, and therefore $J^{m_k}(C) \cap V \neq \emptyset$ for 
large-enough $k$. This contradiction concludes the proof. $\hfill\square$

\vspace{0.4cm}

Now fix $x_0 \in X$, and let $v_0 := ctr(M_{x_0})$.

\vspace{0.02cm}

\begin{lema} \label{dorila}
If $(f_k)$ is a sequence of group elements such that
$f_k (x_0) \rightarrow x_0$, then $I(f_k,x_0) v_0 \rightarrow v_0$.
\end{lema}

\noindent{\bf Proof.} Since $I(f_k,x_0) (M_{x_0})$ converges to $M_{x_0}$, by
Proposition \ref{se-usa}, we have that
$ctr(I(f_k,x_0) (M_{x_0}))$ converges to $ctr(M_{x_0})$.
Since the map \hspace{0.01cm} $ctr$ \hspace{0.01cm} commutes with isometries,
$$I(f_k,x_0) v_0 = I(f_k,x_0) \big( ctr(M_{x_0}) \big) \longrightarrow ctr(M_{x_0}) 
= v_0,$$
thus showing the lemma. $\hfill\square$

\vspace{0.32cm}

Denote the closure of the orbit of $(x_0,v_0)$ by $\hat{M}$. This is a compact 
invariant set. Moreover, by Lemma \ref{dorila}, the fiber $\hat{M}_{x_0}$ is 
reduced to $v_0$.

\vspace{0.1cm}

\begin{lema} \label{cuarterola} For each $x \in X$, the set $\hat{M}_{x}$
is reduced to a single point.
\end{lema}

\noindent{\bf Proof.} Assume for a contradiction that $\hat{M}_{x}$ contains
at least two points, say $v \neq v^*$. Since the $\Gamma$-action on $X$ is minimal,
there exists a sequence $(f_k)$ in $\Gamma$ such that $f_k (x) \rightarrow x_0$. 
The points $(f_k(x), I(f_k,x) v)$ and $(f_k(x), I(f_k,x) v^*)$ belong to 
$\hat{M}$ for each $k$. Passing to a subsequence if necessary, 
we may assume that they converge to $(x_0,w)$ and $(x_0,w^*)$, 
respectively. Notice that because $\hat{M}$ is invariant and closed, 
these two limit points are contained in $\hat{M}$, hence both $w$ 
and $w^*$ are in $\hat{M}_{x_0}$. Now for all $k$, we have
$$d \big( I(f_k,x) v , I(f_k,x) v^* \big) = d(v,v^*).$$
Therefore, $d(w , w^*) = d(v , v^*) > 0$.
However, this contradicts the fact that
$\hat{M}_{x_0} = \{v_0\}$. $\hfill\square$

\vspace{0.55cm}

\noindent{\bf End of the proof.} By Lemma \ref{cuarterola}, the set $\hat{M}$
is the graph of a well-defined function $\varphi: X \rightarrow \mathcal{H}$.
Since $\hat{M}$ is compact, this function is continuous. Finally, because 
$\hat{M}$ is invariant, the curve $x \to (x,\varphi(x))$ satisfies all the 
desired properties. This concludes the proof of the Main Theorem provided 
we give a

\vspace{0.38cm}

\noindent{\bf Proof of Proposition \ref{se-usa}.} 
Given $\varepsilon > 0$, let $B_{\varepsilon}$ be a set within Hausdorff
distance \hspace{0.01cm} $dist_H (B,B_{\varepsilon}) \leq \varepsilon$ \hspace{0.01cm}
from $B$. From the inclusions $B \subset \overline{\B (ctr(B),r_B)}$ and
$B_{\varepsilon} \subset \overline{\B (B,\varepsilon)}$, we obtain
$B_{\varepsilon} \subset \overline{\B (ctr(B),r_B + \varepsilon)}.$
Similarly, we have
$B \subset \overline{\B (ctr(B_{\varepsilon}),r_{B_{\varepsilon}} + \varepsilon)}.$
As a consequence,
\begin{equation}\label{obvio}
\big| r_B - r_{B_{\varepsilon}} \big| \leq \varepsilon.
\end{equation}
Let $m_{\varepsilon}$ be the midpoint between $ctr(B)$ and $ctr(B_{\varepsilon})$.
For each $w \in B_{\varepsilon}$, the median ineequality ({\em i.e.} 
the CAT(0) property) yields
$$d( m_{\varepsilon} , w )^2 = \frac{d( ctr(B_{\varepsilon}) , w )^2}{2} +
\frac{d( ctr(B) , w )^2}{2} - \frac{d( ctr(B_{\varepsilon}) , ctr(B) )^2}{4}.$$
Taking the supremum over all $w \in B_{\varepsilon}$ and using (\ref{obvio}),
we obtain
\begin{eqnarray*}
r^2_{B_{\varepsilon}} 
&=& \sup_{w \in B_{\varepsilon}} d \big( m_{\varepsilon} , w \big)^2\\
&\leq& \frac{r^2_{B_{\varepsilon}}}{2} +
\frac{\sup_{w \in B_{\varepsilon}} d( ctr(B) , w )^2}{2}
- \frac{d( ctr(B_{\varepsilon}) , ctr(B) )^2}{4}\\
&\leq& \frac{r^2_{B_{\varepsilon}}}{2} + \frac{1}{2} \left[ \sup_{w \in B} d( ctr(B) , w )
+ dist_{H} (B,B_{\varepsilon}) \right]^2
- \frac{d( ctr(B_{\varepsilon}) , ctr(B) )^2}{4}\\
&\leq& \frac{r^2_{B_{\varepsilon}}}{2} + \frac{1}{2} \big[ r_{B} + \varepsilon \big]^2
- \frac{d( ctr(B_{\varepsilon}) , ctr(B) )^2}{4},
\end{eqnarray*}
hence
$$d \big( ctr(B_{\varepsilon}) , ctr(B) \big)^2 \hspace{0.1cm} \leq \hspace{0.1cm} 
4 \Big( \frac{\big[ r_{B} + \varepsilon \big]^2}{2} - \frac{r^2_{B_{\varepsilon}}}{2} \Big) 
\hspace{0.1cm} \leq \hspace{0.1cm} 
2 \big( [r_B + \varepsilon]^2 - [r_B - \varepsilon]^2 \big) 
\hspace{0.1cm} = \hspace{0.1cm} 8 \varepsilon r_{B}.$$
Since the right-side expression converges to zero together with 
$\varepsilon$, this concludes the proof. $\hfill\square$


\subsubsection{Second proof}

\hspace{0.45cm} We next provide an even more geometric argument of proof. Consider 
a compact, invariant set $M$ for the skew action of $\Gamma$. Fix $x_0 \in X$, 
and denote $v_0 := ctr (M_{x_0})$. Denote also the closure of the orbit of 
$(x_0,v_0)$ by $\hat{M}$. The main step in the first proof was to show 
that $\hat{M}_{x_0}$ is reduced to $\{ v_0 \}$. (Starting from this, 
Lemma \ref{cuarterola} shows that $\hat{M}_x$ is reduced to a single point 
for each $x \in X$, which allows to conclude in the same way as before.) 

Assume for a contradiction that $\hat{M}_{x_0}$ contains a point  
$v_0'$ distinct from $v_0$, and let $r_0 := r_{{M}_{x_0}}$ and 
$\varepsilon_0 := d(v_0' , v_0 ) > 0$. There must be a sequence 
$(f_k)$ in $\Gamma$ such that $( f_k (x_0), I(f_k,x_0) v_0 )$ 
converges to $(x_0,v_0')$. Since ${M}_{f_k (x_0)} \subset 
\overline{\B ( I(f_k,x_0) v_0, r_0 )}$, given 
$\varepsilon > 0$, we must have, for large-enough $k$, 
\begin{equation}
\label{posicion-1}
{M}_{f_k (x_0)} \subset \overline{\B \big( v_0', r_0 + \varepsilon \big)}.
\end{equation}
We now claim that for large-enough $k$, we also have 
\begin{equation}
\label{posicion-2}
{M}_{f_k (x_0)} \subset \overline{\B (v_0, r_0 + \varepsilon)}.
\end{equation}
Indeed, if not, then there would be a sequence $(v_n)$ such that  
$v_{k_n}$ belongs to ${M}_{f_{k_n} (x_0)} \setminus \B (v_0, r_0 + \varepsilon)$ 
for an increasing sequence of integers $(k_n)$. Passing to a subsequence 
if necessary, this would yield a limit point 
$v^* \in {M}_{x_0} \setminus \B (v_0, r_0 + \varepsilon)$, 
which is absurd. 

Now, (\ref{posicion-1}) and (\ref{posicion-2}) yield 
(for a large-enough $k$ depending on $\varepsilon > 0$)
$${M}_{f_k (x_0)} \subset \overline{\B ( v_0', r_0 + \varepsilon )} 
\bigcap \overline{\B (v_0, r_0 + \varepsilon)}.$$
The contradiction we seek comes from the fact that right-side set has 
radius at most $r_0 - \varepsilon$ provided that $\varepsilon$ is less  
than or equal to $\varepsilon_0^2 / 16 r_0$. Indeed, letting $v$ be the 
midpoint between $v_0$ and $v_0'$, the median inequality yields, for each 
$w \in \overline{\B ( v_0', r_0 + \varepsilon )} 
\bigcap \overline{\B (v_0, r_0 + \varepsilon)}$,
$$d( w , v )^2 + \frac{d( v_0 , v_0' )^2}{4} \leq  
\frac{d( w , v_0 )^2}{2} + \frac{d( w , v_0' )^2}{2}.$$
Thus, 
$$d( w , v )^2 \leq \frac{(r_0 + \varepsilon)^2}{2} + 
\frac{(r_0 + \varepsilon)^2}{2} - \frac{\varepsilon_0^2}{4} \leq (r_0 - \varepsilon)^2.$$ 
Therefore, the set \hspace{0.05cm} $\overline{\B ( v_0', r_0 + \varepsilon )} 
\bigcap \overline{\B (v_0, r_0 + \varepsilon)}$ \hspace{0.05cm} is contained 
in \hspace{0.05cm} $\overline{\B (v,r_0 - \varepsilon)}$, \hspace{0.05cm} 
which shows that its radius is at most $r_0 - \varepsilon$.


\subsubsection{Third proof}

\hspace{0.45cm} This proof is restricted to the case of cocycles of isometries of 
$\mathbb{R}^{\ell}$, but very likely it extends to general proper CAT(0) spaces. 
Its interest relies in that it relates previous discussion to a classical notion.

\vspace{0.25cm}

\noindent{\bf The recurrence semigroup.} Let us consider a general skew 
action of a semigroup $\Gamma$ on $X \times \mathbb{R}^{\ell}$, namely 
$f \!: (x,v) \to (f(x),I(f,x)v)$, so that the $\Gamma$-action on $X$ is 
minimal and each $I(f,x)$ is an isometry of $\mathbb{R}^{\ell}$. Given $x \in X$, 
we denote by $R_{x}$ the set of isometries $I$ of $\mathbb{R}^{\ell}$ such that 
$I = \lim_{k} I(f_k,x)$ for a sequence of elements $f_k \in \Gamma$ 
satisfying $f_k (x) \rightarrow x$. We begin with the following

\vspace{0.1cm}

\begin{lema} \label{fin-semigroup}
Assume that there is a compact subset $K$ of $Isom(\mathbb{R}^{\ell})$ such that 
$I(f,x)$ lies in $K$ for every $f \!\in\! \Gamma$ and all $x \!\in\! X$. (This is 
equivalent to that the set $I(f,x) v$ is bounded for each $v \in \mathbb{R}^{\ell}$.) 
Then for every $x \in X$, the set $R_{x}$ is a semigroup.
\end{lema}

\noindent{\bf Proof.} Let $d$ be the metric on $X$, and let $dist$ be the left-invariant 
distance on the group of isometries of $\mathbb{R}^{\ell}$ induced by \esp 
$dist (\Psi + \rho, Id) = \| \Psi - Id \| + \| \rho \|.$ \esp One readily 
checks that there is a constant $C = C_{K}$ such that $dist$ is perturbed 
under right-translation by a factor at most $C$, that is, \esp 
$dist(I_1 I, I_2 I) \leq C dist (I_1,I_2)$ \esp for all 
$I \in K$ and all $I_1,I_2$ in $Isom(\mathbb{R}^{\ell})$.  

Given $I_1,I_2$ in $R_x$, we need to show that $I_1 I_2$ also 
belongs to $R_x$. For $i \in \{1,2\}$, choose a sequence $(f_{i,k})_{_k}$ such that 
$f_{i,k} (x) \rightarrow x$ and $I(f_{i,k},x) \to I_i$. Given $\varepsilon > 0$, 
there is an integer $k_1 \in \mathbb{N}$ such that, for all $k \geq k_1$,
$$d \big( f_{1,k} (x), x \big) \leq \varepsilon \quad 
\mbox{ and } \quad dist \big( I(f_{1,k},x),I_1 \big) \leq \varepsilon.$$
By continuity, there exists $\delta \!\in ]0,\varepsilon[$ 
such that if $d(x,y) \leq \delta$, then  
$$d \big( f_{1,k_1}(x),f_{1,k_1}(y) \big) \leq \varepsilon \quad \mbox{ and } \quad 
dist \big( I(f_{1,k_1},x), I(f_{1,k_1},y) \big) \leq \varepsilon.$$
Fix $k_2 \in \mathbb{N}$ large-enough that 
$$d \big( f_{2,k_2}(x),x \big) \leq \delta \quad \mbox{ and } \quad 
dist \big( I(f_{2,k_2},x), I_2 \big) \leq \varepsilon.$$
We have 
$$d \big( f_{1,k_1} f_{2,k_2} (x), x \big) \leq 
d \big( f_{1,k_1} (f_{2,k_2}(x)), f_{1,k_1}(x) \big) + d \big( f_{1,k_1}(x),x \big) 
\leq 2 \varepsilon.$$
Moreover, using the almost invariance of \hspace{0.01cm} $dist$, \hspace{0.01cm} from \hspace{0.06cm} 
$dist ( I(f_{1,k_1},x), I(f_{1,k_1},f_{2,k_2}(x)) ) \leq \varepsilon$, 
\hspace{0.06cm} we get 
\begin{small}
$$dist \big( I(f_{1,k_1},x) I(f_{2,k_2},x), I(f_{1,k_1} f_{2,k_2},x) \big) = 
dist \big( I(f_{1,k_1},x) I(f_{2,k_2},x), I(f_{1,k_1},f_{2,k_2}(x)) I(f_{2,k_2},x) \big) 
\leq C \varepsilon.$$
\end{small}Therefore, 
\begin{small}
\begin{eqnarray*}
dist \big( I_1 I_2, I(f_{1,k_1} f_{2,k_2},x) \big) 
\! &\leq& \! dist \big( I_1 I_2, I(f_{1,k_1},x) I(f_{2,k_2},x) \big) + 
dist \big( I(f_{1,k_1},x) I(f_{2,k_2},x), I(f_{1,k_1} f_{2,k_2},x) \big)\\
\! &\leq& \! dist \big( I_1 I_2, I(f_{1,k_1},x) I_2 \big) + 
dist \big( I(f_{1,k_1},x) I_2, I(f_{1,k_1},x) I(f_{2,k_2},x) \big) + C \varepsilon\\
\! &=& \! C dist \big( I_1, I(f_{1,k_1},x) \big) + 
dist \big( I_2, I(f_{2,k_2},x) \big) + C \varepsilon\\
\! &\leq& \! (2 + C) \varepsilon.
\end{eqnarray*}\end{small}Summarizing, for each $\varepsilon > 0$, we have 
found an element $f \in \Gamma$, namely, $f := f_{1,k_1} f_{2,k_2}$, such that
$$d \big( f(x),x \big) \leq 2 \varepsilon \quad \mbox{ and } \quad 
dist \big( I_1 I_2, I(f,x) \big) \leq (2 + C) \varepsilon.$$
By definition, this shows that $I_1 I_2$ belongs to $R_x$. $\hfill\square$

\vspace{0.35cm}

We will call $R_x$ the {\em recurrence semigroup of} $x$. (A closely 
related notion was developed in \cite{atkinson}.) We must emphasize that, 
in general, $R_x$ is not a group, even if $\Gamma$ has a group structure. 
Nevertheless, if $\Gamma$ is a group and its action on the basis $X$ is 
{\em equicontinuous}, then $R_x$ is a group. Indeed, given $I \in R_x$, choose 
a sequence $(f_k)$ in $\Gamma$ so that $f_k (x) \to x$ and $I(f_k,x) \to I$. 
By equicontinuity, we also have $f_k^{-1}(x) \to x$, and passing to a subsequence 
if necessary, we may assume that $I(f^{-1}_k,x)$ converges to an isometry $I'$. 
Now, the proof of Lemma \ref{fin-semigroup} yields that $I(f_n^{-1} f_m) \to I'I$ 
as $n,m$ go to infinity. Letting $n=m$ go to infinity, this obviously implies  
$Id = I'I$, that is, $I' \in R_x$ is the inverse of $I$.

\vspace{0.1cm}

\begin{lema} \label{cle-cle} If $M$ is the closure of the orbit of a point 
$(x,v) \in X \times \mathbb{R}^{\ell}$ under the skew action, then the set $M_x := 
\{w \in \mathbb{R}^{\ell} \!: (x,w) \in M\}$ coincides with $\{ I v \!: I \in R_x \}$. 
\end{lema}

\noindent{\bf Proof.} Each point in $M$ is of the form $\lim_k ( f_k(x),I(f_k,x)v )$ 
for a sequence of elements $f_k \in \Gamma$. Thus, each point of $w \in M_x$ has 
the form $\lim_k I(f_k,x)v$ for a sequence $(f_k)$ in $\Gamma$ such that 
$f_k (x) \rightarrow x$. The lemma follows from the fact that the set of isometries 
sending a prescribed vector $v \in \mathbb{R}^{\ell}$ into some fixed bounded neighborhood 
of another prescribed vector $w \in \mathbb{R}^{\ell}$ is compact. $\hfill\square$

\vspace{0.38cm}

Assume now that the orbit of some point $(x,v) \in X \times \mathbb{R}^{\ell}$ is bounded, 
and let $M$ be its closure. Fix a point $x_0 \in X$, and let $v_0:= ctr(M_{x_0})$. 
Finally, let $\hat{M}$ be the closure of the orbit of the point $(x_0,v_0)$.

\vspace{0.1cm}

\begin{lema} \label{pingana} The set 
$\hat{M}_{x_0} \! := \! \{ w \in \mathbb{R}^{\ell} \! : (x_0,w) \in \hat{M} \}$ 
reduces to $\{v_0\}$.
\end{lema}

\noindent{\bf Proof.} By Lemma \ref{cle-cle}, the point $v_0$ 
may be written as $ctr ( \{Iv \!: I \in R_{x_0} \})$. By the 
semigroup version of the Bruhat-Tits center lemma, this point 
is fixed by every element of $R_{x_0}$. In other words, the set 
$\{I v_0 \!: I \in R_{x_0} \}$ reduces to $\{v_0\}$. Finally, by Lemma 
\ref{cle-cle} again, this set coincides with $\hat{M}_{x_0}$. $\hfill\square$

\vspace{0.4cm}

The rest of the third proof works as the final part of the first one. Indeed, as in 
Lemma \ref{cuarterola}, one may show that for each $x \in X$, the set $\hat{M}_{x}$ 
is reduced to a single point. Hence, the set $\hat{M}$ is the graph of a 
well-defined function $\varphi: X \rightarrow \mathbb{R}^{\ell}$, and the 
curve $x \mapsto (x,\varphi(x))$ satisfies all the desired properties. 


\subsubsection{Fourth proof}

\hspace{0.45cm} Assume once again that there is a nonempty, compact, 
forward-invariant set $M$ for the skew action of $\Gamma$. 

\begin{lema} \label{radio-cte}
The function $D \!: X \to [0,\infty[$ that makes correspond, to each 
$x \in X$, the diameter of the set $M_x := \{w \in \mathcal{H} \!: (x,w) \in M \}$, 
is constant.
\end{lema}

\noindent{\bf Proof.} The function $D$ is invariant under the $\Gamma$-action on $X$. 
Since this action is assumed to be minimal, in order to prove that $D$ is constant, 
it suffices to show that it is upper-semicontinuous. To do this, let $(x_n)$ be an 
arbitrary sequence of points converging to a certain $x \in X$. Let $(x_{n_k})$ be 
a subsequence such that $\lim_{k} D(x_{n_k}) = \lim\sup_n D(x_n) =: \Delta$. For 
each $k$, let $(v_k,w_k)$ be a pair of points of $M_{x_{_{n_k}}}$ at distance 
$D(x_{n_k})$. Passing to a subsequence if necessary, we may suppose that $v_k$ 
(resp. $w_k$) converges to a certain $v \in M_x$ (resp. $w \in M_x$). Clearly,  
$d( v , w ) \geq \lim_k d(v_k , w_k ) = \Delta$. In particular, the diameter of 
$M_x$ is greater than or equal to $\Delta$. This shows the lemma. $\hfill\square$

\vspace{0.5cm}

We will denote by $D(M)$ the common value of the diameter of the fibers $M_x$. Notice 
that, denoting by $cv(M)$ the {\em convex closure} of $M$ along the fibers, we have 
$D \big( cv(M) \big) = D(M)$. Moreover, straightforward arguments show that $cv(M)$ 
is also compact and invariant. 

\vspace{0.1cm}

Now consider the set 
\begin{small}
$$M^* := \big\{(x,v) \!: \mbox{ there exist } v_1,v_2 \mbox{ at distance } D \mbox{ in } M_x 
\mbox{ such that } v \mbox{ is the midpoint of the segment } \overline{v_1v_2} \big\}.$$ 
\end{small}Notice that $M^* = cv(M)^*$ and $cv(M^*)$ is contained in $cv(M)$. Moreover, $M^*$ 
is invariant under the skew action. Furthermore, easy compactness-type arguments show that  
$M^*$ is closed (hence compact) and nonempty. ({\em A priori}, the fibers of $M^*$ do not 
vary continuously.) Finally, the preceding lemma applied to $M^*$ shows that all its 
fibers have the same diameter. 

The next lemma is a direct consequence of \cite[Lemma 3.2.3]{BT}, and we reproduce 
the proof just for the reader's convenience.

\vspace{0.1cm}

\begin{lema} \label{diameter}
One has the inequality $D (M^*) \leq D(M) / \sqrt{2}$. Moreover, this estimate is sharp.
\end{lema}

\noindent{\bf Proof.} Fix $x \in X$ 
and let $v,w$ be points in $M^*_x$. By definition, there exist two pairs 
of points $(v_1,v_2)$ and $(w_1,w_2)$ at distance $D(M)$ in $M_x$ such that $v$ (resp. $w$) 
is the midpoint of $\overline{v_1v_2}$ (resp. $\overline{w_1w_2}$). The median inequality 
applied to the triangle $\Delta(v_1,w_1,w_2)$ yields 
$$d(v_1 , w)^2 + \frac{D(M)^2}{4} \leq  
\frac{d(v_1 , w_1)^2}{2} + \frac{d(v_1 , w_2)^2}{2} 
\leq \frac{D(M)^2}{2} + \frac{D(M)^2}{2} = D(M)^2,$$
hence, \hspace{0.01cm} $d(v_1 , w )^2 \leq 3 D(M)^2 / 4$. \hspace{0.01cm} 
Similarly, \hspace{0.01cm} $d(v_2 , w )^2 \leq 3 D(M)^2 / 4$. \hspace{0.01cm} 
Using this, the median inequality for $\Delta(v_1,v_2,w)$ yields 
$$d(v , w)^2 + \frac{D(M)^2}{4} \leq \frac{d( v_1 , w )^2}{2} + 
\frac{d( v_2 , w )^2}{2} \leq \frac{3D(M)^2}{4}.$$
This easily leads to the estimate of the lemma. To see that this estimate is sharp, 
it suffices to consider the case where each $M_x$ consists of four points that are 
the vertices of a tetrahedron. (In dimension 2, the constant $\sqrt{2}$ can be 
replaced by $2$.) $\hfill\square$

\vspace{0.5cm}

\noindent{\bf End of the proof.} Let us define the sequence of nonempty, compact, invariant 
sets $M_n$ by $M_1 := cv(M)$ and $M_n : = cv(M_{n-1}^*)$ for each $n \geq 2$. Since 
$M_n \subset M_{n-1}$ holds for each $n > 1$, the set  $\hat{M} := \bigcap_{n \geq 1} M_n$ 
is also nonempty and compact, as well as invariant. Moreover, Lemma \ref{diameter} 
implies that $D(\hat{M}) = 0$. In other words, each fiber $\hat{M}_x$ consists of 
a single point $\varphi(x)$, and the thus-defined function $\varphi$ satisfies all 
the desired properties. (Its continuity follows from the fact that its graph is compact.)

\vspace{0.2cm}

\begin{obs} {\em It is very instructive to compare the technique of the preceding proof 
with the three previous ones. Given a bounded subset $B \subset \mathcal{H}$, we let 
$B_1 := B$, and having defined $B_2,\ldots,B_{n-1}$, we let $B_n$ be the set of midpoints 
of segments between points of $B_{n-1}$ situated at distance $diam (B_{n-1})$. Finally, 
we call the point $ctr^*(B) := \bigcap _{n \geq 1} B_n$ the Bruhat-Tits center of $B$. 

In general, $ctr^*(B)$ does not coincide with $ctr(B)$. For example, if $B$ consists of three 
points that are the vertices of a triangle $\Delta$ all of whose angles are $\leq \pi$, then 
$ctr (B)$ coincides with the circumcenter of $\Delta$. However, if the sides of $\Delta$ 
have different length, then $ctr^*(B)$ is nothing but the midpoint of the largest side.}
\end{obs}


\subsection{The case of infinite-dimensional Hilbert space fibers}
\label{caso-infinito}

\subsubsection{The lack of continuity of invariant sections given by centers 
along the fibers}
\label{ejem-disc}

\hspace{0.45cm} For fibers that are infinite dimensional 
Hilbert spaces, none of the strategies of proof proposed so far works. 
Indeed, there is a serious technical problem in defining the recurrence semigroup  
(the group of linear isometries is not compact when endowed with the norm-topology). 
Moreover, the continuity in the Hausdorff topology for the weak topology does not 
guarantee the continuity of the center. Finally, the diameter of the fibers of an 
skew-invariant, weakly-compact set is not necessarily constant. 

In a more concrete way, the example below showing that the center along the fibers of a 
weakly-compact invariant set may fail to be continuous illustrates all these technical 
problems. 

\begin{ejem} \label{shifteado}
{\em Let us consider the Hilbert space $\mathcal{H} \sim \ell^{2} (\mathbb{Z})$, 
and let $\Gamma \sim \mathbb{Z}$ be acting on $X$ by powers of a minimal 
homeomorphism $T$. Let us consider the skew action with linear part induced by 
$\Psi(1,x) v_n = v_{n+1}$ for every $x$, where $\{v_n\}$ is an orthonormal basis 
of $\mathcal{H}$, and with translation part $\rho \!: X \to \mathcal{H}$ vanishing  
everywhere. Fix $x_0 \in X$, and consider the two-points set $\{ (x_0,0), (x_0,v_0) \}$. 
The closure (for the weak topology) of its orbit under the skew action is a set $M$ whose 
fiber over $x \in X$ coincides with $\{ 0 \}$ if $x$ is not in the orbit of $x_0$, and with 
$\{0, v_n \}$ if $x = T^n (x_0)$. In the first case, we have $ctr (M_x) = \{ 0 \}$, 
while for $x = T^n (x_0)$ we have $ctr (M_x) = \{ v_n / 2\}$. It is then easy to 
see that the function $\varphi \!: x \mapsto ctr (M_x)$ is not weakly continuous 
(namely, it is discontinuous at every point).} 
\end{ejem}


\subsubsection{Existence of weakly-continuous invariant sections}
\label{weak-cont}

\hspace{0.45cm} In this section, we deal with a skew action by isometries of a 
semigroup $\Gamma$ so that the fibers are a Hilbert space $\mathcal{H}$ and the 
dynamics on the basis is minimal. We assume 
that, for all $f \in \Gamma$, the map $I(f,\cdot): X \rightarrow Isom (\mathcal{H})$ is 
continuous for the strong topology. Writing $I(f,\cdot) = \Psi(f,\cdot) + \rho(f,\cdot)$, 
this means that $\Psi(f,\cdot): X \rightarrow U (\mathcal{H})$ is norm-continuous and 
$\rho(f,\cdot): X \rightarrow \mathcal{H}$ is continuous for the strong topology on 
$\mathcal{H}$. We endow $X \times \mathcal{H}$ with the product topology, where the 
topology on the factor $\mathcal{H}$ is the weak one.

Suppose that there exists a bounded orbit for the skew action, and let us 
consider its convex closure $M$. By this we mean the smallest compact set 
that contains the given set and is convex along the fibers, in the sense that 
if $(x,v), (x,w)$ belong to $M$ then $( x, \lambda v + (1-\lambda)w ) \in M$ for all 
$0 \leq \lambda \leq 1.$ The family $\mathcal{F}$ of nonempty, compact, invariant sets 
that are convex along the fibers is ordered by inclusion. A straightforward application 
of Zorn's lemma shows that it contains a minimal element. The crucial step is the next 

\vspace{0.1cm}

\begin{lema} \label{crucial} For each minimal element $M$ of the family $\mathcal{F}$, 
the fiber $M_x$ above $x$ consists of a single vector, for each $x \in X$.
\end{lema}

\vspace{0.1cm}

This lemma yields a weakly-continuous invariant section for the skew action. Indeed, 
the set $M$ will be the graph of a function $\varphi \!: X \to \mathcal{H}$ which is 
weakly continuous, since its graph is compact. Moreover,  since the set $M$ is 
invariant, $\varphi$ is an invariant function.

\vspace{0.45cm}

\noindent{\bf Proof of Lemma \ref{crucial}.} Assume that for some $x_0 \in M$ the fiber 
$M_{x_0}$ contains two vectors $v_1,v_2$ at a distance $\| v_1-v_2 \| =: \varepsilon > 0$. 
Let $r(M) > 0$ be the infimum of the radius $r$ such that, for all $x \in X$, the fiber 
$M_x$ is contained in $\mathrm{Ball}(0,r)$. Given $\kappa < 1$, there must exist 
$(y,w) \in M$ such that \esp $\| w \| \geq \kappa r(M).$ \esp Fix such a 
$w \in M_{y}$ and $\kappa < 1$ such that 
$$\kappa > \sqrt[4]{1 - \frac{\varepsilon^2}{4 r(M)^2}}.$$
For each $\lambda < 1$, let $P_{\lambda}$ be the affine hyperplan 
$\lambda w + \langle w \rangle^{\perp}$. This hyperplane divides 
the whole fiber $\mathcal{H}$ above $y$ into two closed 
hemispheres $P^+_{\lambda}, P_{\lambda}^-$, where 
$w$ belongs to the interior of $P^+_{\lambda}$. 

Let $u$ be the midpoint between $v_1$ and $v_2$. By convexity, the point 
$(x_0,u) := (x, (v_1+v_2)/2)$ must belong to $M$. We claim that the closure 
of its orbit must intersect the hemisphere $\{ y \} \times P^+_{\lambda}$. 
Otherwise, the convex closure of its orbit would be a nonempty, compact, 
invariant set that is convex along the fibers and it is strictly contained 
in $M$ (it does not contain $(y,w)$). However, this contradicts the fact 
that $M$ is a minimal element of $\mathcal{F}$.

We thus conclude that for each $\lambda_* < \lambda$ there exists $f \in \Gamma$ 
such that $I(f,x_0) u \in P_{\lambda_*}^+$. Fixing such a $\lambda_*$ so that 
$$\lambda_* > \sqrt[4]{1 - \frac{\varepsilon^2}{4 r(M)^2}},$$
we claim that 
\begin{equation}\label{se-sale}
\mbox{either } \, I(f,x)v_1 \, \mbox{ or } \, I(f,x)v_2 \,
\mbox{ lies outside } \overline{\mathrm{Ball} \big( 0,r(M) \big)}.
\end{equation}
Before proving this claim, notice that it contradicts the definition of $r(M)$, 
thus concluding the proof.

The proof of (\ref{se-sale}) relies on the uniform convexity of $\mathcal{H}$. 
In a quantitative manner, since $I(f,x) u$ lies in $P_{\lambda_*}^+$, its 
norm is at least $\kappa \lambda_* r(M)$. By the median equality 
$$\frac{\big\| I(f,x)v_1 - I(f,x)v_2 \big\|^2}{4} 
= \frac{\big\| I(f,x) v_1 \big\|^2}{2} 
+ \frac{\big\| I(f,x) v_2 \big\|^2}{2} - \big\| I(f,x)u \big\|^2,$$
this yields 
$$\frac{\varepsilon^2}{4} \leq \frac{\big\| I(f,x) v_1 \big\|^2}{2} 
+ \frac{\big\| I(f,x) v_2 \big\|^2}{2} - \kappa^2 \lambda_*^2 r(M)^2.$$
Assuming that both $I(f,x)v_1$ and $I(f,x)v_2$ are in 
$\overline{\mathrm{Ball}(0,r(M))}$, this implies that  
$$\frac{\varepsilon^2}{4} \leq r(M)^2 - \kappa^2 \lambda_*^2 r(M)^2 = 
r(M)^2 \big( 1 - \kappa^2 \lambda_*^2 \big),$$
which contradicts our choice of the constants $\kappa,\lambda_*$. $\hfill\square$
  

\subsubsection{Strong continuity of weakly-continuous invariant sections}
\label{weak/strong}

\hspace{0.45cm} In order to complete the proof of the Main Theorem in the infinite-dimensional 
case, we need to show that in the context of \S \ref{weak-cont}, the following hods:

\vspace{0.1cm}

\begin{prop} \label{prop:weak-implies-strong}
Every weakly-continuous solution of the cohomological equation {\em (\ref{chupalo-vos})} 
is strongly continuous. Equivalently, all weakly-continuous, skew-invariant sections are 
strongly continuous.
\end{prop}

\vspace{0.1cm}

To show this proposition, we begin by giving a geometrical criterion for strong 
continuity. To do this, given a function $\varphi : X \rightarrow \mathcal{H}$, 
we define its {\em oscillation} at a point $x$ as 
$$\mathrm{osc} (\varphi)(x) := \limsup_{\{y,z\} \to \{x\}} 
\big\| \varphi(y) - \varphi(z) \big\|.$$
Our first lemma should be clear from the definition.

\vspace{0.1cm}

\begin{lema} \label{un} 
The map $\varphi$ is strongly continuous at a point 
$x \in X$ if and only if $\osc (\varphi) (x) = 0$.
\end{lema}

\vspace{0.3cm}

Our second lemma involves the underlying dynamics of our setting.

\vspace{0.1cm}

\begin{lema} \label{trois} If the curve $x \mapsto ( x,\varphi(x) )$ 
is skew invariant, then the function $x \mapsto \osc (\varphi)(x)$ 
is invariant under the action of $\Gamma$ on $X$. 
\end{lema}

\noindent{\bf Proof.} Let $f \in \Gamma$  and $x \in X$ be given. It is enough 
to show that, given sequences $(y_n)$, $(z_n)$ converging to $x$ so that 
$\| \varphi(y_n) - \varphi(z_n) \|$ converges to some value $\varepsilon$, 
then there exist $(\bar{y}_n)$ and $(\bar{z}_n)$ converging to $f(x)$ so that 
$\| \varphi(\bar{y}_n) - \varphi(\bar{z}_n) \|$ also converges to $\varepsilon$. 
We will show that this holds for $\bar{y}_n := f (y_n)$ and $\bar{z}_n := f (z_n)$. 
Indeed, the value of 
$$\big\| \varphi \big( f(y_n) \big) - \varphi \big( f(z_n) \big) \big\|$$
may be written as
$$\big\| I (f,y_n) \varphi(y_n) - I (f,z_n) \varphi(z_n) \big\|,$$ 
and differs from 
$$\big\| I(f,x) \varphi(y_n) - I(f,x) \varphi(z_n) \big\| = \| y_n - z_n \|
$$
by no more than 
$$\big\| I(f,y_n) \varphi(y_n) - I (f,x) \varphi(y_n) \big\| + 
\big\| I(f,z_n) \varphi(z_n) - I (f,x) \varphi(z_n) \big\|,$$ 
which equals 
$$\big\| \big( \Psi(f,y_n) - \Psi(f,x) \big) \varphi(y_n) + \rho(y_n) - \rho(x) \big\| +
\big\| \big( \Psi(f,z_n) - \Psi(f,x) \big) \varphi(z_n) + \rho(z_n) - \rho(x) \big\|.$$
Since $\varphi$ is weakly continuous, it must be bounded, say by a constant 
$C > 0$. This implies that the last expression above is bounded from above by 
$$C \big( \| \Psi(f,y_n) - \Psi(f,x) \| + \| \Psi(f,z_n) - \Psi(f,x) \| \big) 
+ \| \rho(y_n) - \rho (x)\| + \| \rho(z_n) - \rho(x) \|.$$
By the norm-continuity of $\Psi(f,\cdot)$ and the strong continuity of 
$\rho(f,\cdot)$, this last expression converges to zero. 
This concludes the proof. $\hfill\square$ 

\vspace{0.5cm}

The next lemma shows that the function 
$x \mapsto \osc (\varphi)(x)$ is lower-semicontinuous.

\vspace{0.1cm}

\begin{lema} \label{deux} 
For each $\varepsilon > 0$, the set 
$\{x \!: \osc (\varphi)(x) < \varepsilon\}$ is open in $X$.
\end{lema}

\noindent{\bf Proof.} Given $x_0$ in this set, let 
$\varepsilon_0 := \osc (\varphi)(x_0) < \varepsilon$. 
Then there exists $\delta > 0$ such that, for all $y,z$ at distance 
$< \delta$ from $x_0$, we have $\|\varphi(y) - \varphi(z) \| \leq \frac{1}{2} 
(\varepsilon + \varepsilon_0)$. This clearly implies that, for all $x \in X$ 
such that $dist (x,x_0) < \delta$, we have $\osc (\varphi)(x) \leq 
\frac{1}{2} (\varepsilon + \varepsilon_0) < \varepsilon$. In other 
words, the $\delta$-neighborhood of $x_0$ is contained in 
$\{x: \osc (\varphi)(x) < \varepsilon\}$, thus showing the 
lemma. $\hfill\square$

\vspace{0.35cm}

We are now ready to prove Proposition \ref{prop:weak-implies-strong}. Indeed, since 
the $\Gamma$-action on $X$ is assumed to be minimal, Lemmata \ref{trois} and \ref{deux} 
imply that each set $\{x: \osc (\varphi)(x) < \varepsilon\}$ is either empty or 
coincides with the whole space $X$. If we are able to detect a point where 
$\varphi$ is strongly continuous, then 
by Lemma \ref{un} we will have a point in each of these sets. Hence, each of these sets 
will coincide with $X$, so that the oscillation of $\varphi$ at every point will 
be zero. By Lemma \ref{un} again, this will imply that $\varphi$ is strongly continuous.

Thus, to conclude the proof, we need to ensure the existence of a point of 
strong continuity for $\varphi$. This follows from the following well-known 

\vspace{0.18cm}

\begin{lema} Every weakly-continuous function $\varphi: X \rightarrow \mathcal{H}$ 
is strongly continuous on a $G_{\delta}$-set. 
\end{lema}

\noindent{\bf Proof.} Let $\mathcal{H}_1 \subset \mathcal{H}_2 \subset \ldots$ be a 
sequence of finite-dimensional subspaces such that $\mathcal{H}_* := \bigcup_n \mathcal{H}_n$ 
is dense in $\mathcal{H}$. Since $\varphi$ is weakly continuous, each of the functions 
$x \mapsto \| \varphi (x) \|_{\mathcal{H}_n}$ is continuous. By a classical theorem of 
R. Baire \cite{baire} (see also \cite{OU}), 
the pointwise limit of these functions is continuous on a $G_{\delta}$-set 
$X_{\varphi}$. But this pointwise limit is nothing but the function 
$x \mapsto \| \varphi(x) \|$. Recalling now that, in any Hilbert space, 
weak convergence plus convergence of the norm imply strong convergence,  
this yields the strong continuity of $\varphi$ 
on $X_{\varphi}$. $\hfill\square$


\appendix
\section{Appendix. Measurable versus continuous solutions}
\label{measurable-continuous} 

\hspace{0.45cm} We next give a rigidity result for measurable solutions of the 
cohomological equation (\ref{chupalo-vos}) that corresponds to a dynamical 
version/extension of the Corollary to Theorem C. Given a probability measure $\mu$ 
on $X$ that is quasi-invariant under the $\Gamma$-action, we will say that the 
linear part of a skew action on $X \times \mathcal{H}$ is {\em weakly ergodic} 
if the only measurable functions $\phi \!: X \rightarrow \mathcal{H}$ such that 
$\Psi(f,x) \phi(x) = \phi(f(x))$ for all $f \in \Gamma$ and $\mu$-a.e. $x \in X$ 
are the constant ones. 

\begin{ejem} {\em If $\Psi(f,x) = Id$ for all $(f,x)$, then the linear part is 
weakly ergodic if and only if the $\Gamma$-action on $X$ is ergodic w.r.t. $\mu$.}
\end{ejem}

\begin{ejem} \label{rotation-eq} 
{\em In Example \ref{ejem-melnikov}, assume that $T$ is the rotation of angle 
$\alpha \notin \mathbb{Q}$ on the circle (endowed with the Lebesgue measure). 
If $\phi: \mathrm{S}^1 \rightarrow \mathbb{C}$ satisfies 
$\phi(\theta + \alpha) = e^{i \beta} \phi(\theta)$ for a.e. $\theta \in \clo$, 
then $\phi$ must be constant unless $\alpha$ and $\beta$ are rationnaly dependent. 
(See the final argument in Example \ref{ejem-ciclon}.) We thus conclude that the 
linear part of the skew action is weakly ergodic provided $\alpha$ and $\beta$ 
are independent over the rationals.}
\end{ejem}

\begin{ejem} {\em As in Example \ref{shifteado}, 
assume that $\Gamma \sim \mathbb{Z}$ acts by powers of a minimal 
homeomorphism $T$ and that the linear part of its skew action 
on an infinite-dimensional Hilbert space $\mathcal{H}$ is 
generated by $\Psi(1,x) (v_n) = \Psi(v_n) = v_{n+1}$, where 
$\{ v_n \}$ is an orthonormal basis of $\mathcal{H}$ and $x \in X$ is arbitrary. We 
claim that the weak ergodicity holds for any $T$-invariant probability measure $\mu$. 
Indeed, let $\phi(x) = \sum_{n \in \mathbb{Z}} \phi_n(x) v_n$ be a 
measurable function from $X$ to 
$\mathcal{H}$ such that $\Psi (\phi(x)) = \phi(T(x))$, for all $x \in X$. Then we 
have $\phi_{n+1} (T(x)) = \phi_{n} (x)$, for all $x \in X$. If $\phi$ is not 
$\mu$-a.e. equal to zero, then for some $j \in \mathbb{Z}$ and 
$\delta > 0$ we have $\mu(C_{j,\delta}) > 0$,  
where $C_{j,\delta} := \{x \in X : |\phi_j(x)| \geq \delta\}$. By Poincar\'e's recurrence 
theorem, for $\mu$-a.e. point $x \in C_{j,\delta}$ there exists an increasing infinite 
sequence $(n_i)$ such that $T^{-n_i}(x) \in C_{j,\delta}$, hence 
$|\phi_{j+n_i}(x)| = |\phi_j (T^{-n_i}(x))| \geq \delta$. 
However, this is impossible, as $\phi(x)$ belongs to $\mathcal{H}$ 
for $\mu$-a.e. $x \in X$.}
\end{ejem}

\begin{ejem} {\em Let $\Gamma$ be a countable group provided with a spread-out, 
non-degenerate probability distribution $p$, and let $X := P(\Gamma,p)$ be the 
associate Poisson boundary endowed with the corresponding stationary measure 
$\mu$. As a direct consequence of Kaimanovich's double ergodicity theorem \cite{kaimanovich}, 
the linear part of every skew action by isometries of a Hilbert space above 
the natural action of $\Gamma$ on $X$ is weakly ergodic (we assume that $X$ 
is metrizable and compact to fit in our general framework).}
\end{ejem}

\vspace{0.01cm}

\begin{lema} \label{guena} Given a skew action on $X \times \mathcal{H}$ 
whose linear part is weakly ergodic w.r.t. $\mu$, for any two skew-invariant 
{\em measurable} sections $x \mapsto (x,\varphi(x))$ and 
$x \mapsto (x,\overline{\varphi}(x))$, the difference 
$\varphi - \overline{\varphi}$ is a $\mu$-a.e. 
constant vector. If, moreover, there is no common nonzero eigenvector for all 
the $\Psi(f,x)$, then there is at most one skew-invariant measurable solution 
of {\em (\ref{chupalo-vos})}. 
\end{lema}

\noindent{\bf Proof.} For all $f \in \Gamma$, one has $\mu$-a.e. 
$$\varphi \big( f(x) \big) - \Psi(f,x) \varphi(x) = \rho(x) = 
\overline{\varphi} \big( f(x) \big) - \Psi(f,x) \overline{\varphi}(x),$$
hence, 
$$\varphi \big( f(x) \big) - \overline{\varphi} \big( f(x) \big) 
= \big[ \Psi(f,x) \varphi(x) + \rho(x) \big] 
- \big[ \Psi(f,x) \overline{\varphi}(x) + \rho(x) \big] = 
\Psi(f,x) \big( \varphi(x) - \overline{\varphi}(x) \big).$$
Since the linear part of the skew action is assumed to be weakly ergodic, 
this implies that $\varphi - \overline{\varphi}$ is constant. In particular, 
$\varphi - \overline{\varphi}$ is a common eigenvector of all the 
$\Psi(f,x)$. $\hfill\square$ 

\vspace{0.2cm}

\begin{prop} \label{guena-dos} Given a skew action on $X \times \mathcal{H}$, 
assume that the cohomological equation {\em (\ref{chupalo-vos})} admits a 
solution $\varphi \in \mathcal{L}^{\infty}_{\mu}(X,\mathcal{H})$, where 
$\mu$ is such that the linear part is weakly ergodic. If the underlying 
semigroup $\Gamma$ admits a topology with a countable, dense subset such 
that the skew action is continuous, then $\varphi$ is continuous.  
\end{prop}

\noindent{\bf Proof.} For each $f$ lying in a 
countable, dense subset $\Gamma_0$ of $\Gamma$, let \esp 
$Y_f := \big\{ x \in X \!: I(f,x) \varphi(x) \neq \varphi \big( f(x) \big) \big\}.$ 
\esp Then $Y_f$ has null $\mu$-measure, as well as $Y := \bigcup_{f \in \Gamma_0} Y_f$. 
Let $C$ be the essential supremum of the function $x \mapsto \| \varphi(x) \|$, 
and let $Z_0$ be the preimage of $]C,\infty[$ under this function. 
Then $Z_0$ has null $\mu$-measure, as well as 
$Z := \bigcup_{f \in \Gamma_0} f^{-1} (Z_0)$. 
Now, for each $x_0$ in the $\mu$-full measure set 
$X \setminus (Y \cup Z)$ and all $f \in \Gamma_0$, we have 
$$I (f,x_0) \varphi(x_0) = \varphi \big( f(x_0) \big) \quad \mbox{ and } 
\quad \big\| \varphi \big( f(x_0) \big) \big\| \leq C.$$
Since $\Gamma_0$ is dense in $\Gamma$, 
this actually holds for all $f \in \Gamma$, by continuity. 
In other words, the $\Gamma$-orbit of the point 
$(x_0, v_0) := (x_0,\varphi(x_0))$ remains in a bounded subset 
of $\mathcal{H}$. The proposition then follows from 
the Main Theorem combined with Lemma \ref{guena}. $\hfill\square$


\section{Appendix. Cocycles over a shift}
\label{section-shift}

\hspace{0.45cm} Given a minimal homeomorphism $T\!: X \to X$, we consider  
the cocycle of isometries of a Hilbert space $\mathcal{H}$ induced by 
\hspace{0.01cm} $I(1,x) v = \Psi(x) v 
+ \rho(x),$ \hspace{0.01cm} where $\rho: X \to \mathcal{H}$ and 
$\Psi: X \to O(\mathcal{H})$ are continuous. To simplify the notation, 
we write $Tx$ instead of $T(x)$. For each $x \in X$ and 
$k \in \mathbb{N}$, define the $k^{th}$ twisted 
Birkhoff-sum of the cocycle $\rho$ as
$$S_k (\rho)(x) := 
\sum_{i=0}^{k-1} \Psi(T^{k-1}x) \Psi(T^{k-2}x) \cdots \Psi(T^{i+1}x) \rho(T^i x).$$
As it is easy to check, the $k^{th}$ iterate of $(x,v) \in X \times \mathcal{H}$ 
under the skew map $(x,v) \mapsto (Tx, \Psi(x) v + \rho(x))$ coincides with 
$$\big( T^k x, I(k,x)v \big) 
= \big( T^k x, \Psi(T^{k-1}x) \cdots  \Psi(Tx) v + S_k (\rho)(x) \big).$$ 
Assume that there is a bounded orbit for this map, hence a 
continuous solution $\varphi$ to the cohomological equation 
\begin{equation}\label{para-shift}
\varphi(Tx) = \Psi(x) \varphi(x) + \rho(x).
\end{equation}
Then we have
\begin{eqnarray*}
S_k (\rho)(x) 
&=& \sum_{i=0}^{k-1} \Psi(T^{k-1}x) \cdots  \Psi(T^{i+1}x) \rho(T^i x)\\ 
&=& \sum_{i=0}^{k-1} \Psi(T^{k-1}x) \cdots  \Psi(T^{i+1}x)  \big[ \varphi(T^{i+1}x) 
    - \Psi (T^{i} x) \varphi(T^i x) \big]\\
&=& \sum_{i=0}^{k-1} \Psi(T^{k-1}x) \cdots  \Psi(T^{i+1}x) \varphi(T^{i+1}x) - 
    \sum_{i=0}^{k-1} \Psi(T^{k-1}x) \cdots  \Psi(T^{i}x) \varphi(T^{i}x)\\  
&=& \varphi (T^k x) - \Psi(T^{k-1}x) \cdots \Psi(x) \varphi(x).
\end{eqnarray*}
This implies that, as expected, the sequence of functions $S_k (\rho)$ is uniformly 
bounded, hence the orbit of every $(x,v)$ is bounded. Notice that  
if $\Psi$ is constant, then the preceding relation becomes
$$S_k (\rho)(x) = \varphi(T^k x) - \Psi^{k-1} \varphi(x).$$

Let us concentrate on the particular case where 
$\mathcal{H} = \ell^2(\mathbb{Z})$ and $\Psi$ is the (bilateral) 
shift on the canonical basis $\{v_n\}$ of $\mathcal{H}$ (see Example 
\ref{shifteado}). Assuming that $\varphi \!: X \to \mathcal{H}$ solves 
(\ref{para-shift}), fix $x_0 \in X$ and set $x := T^{-k} x_0$ and 
$v := \varphi(T^{-k} x)$. Then we have 
$$\varphi (x_0) = I ( k,T^{-k} x_0 ) \varphi (T^{-k} x_0).$$
Taking the inner product against $v_n$, we get 
\begin{eqnarray*}
\big\langle \varphi (x_0), v_n \big\rangle 
&=& \big\langle \Psi^k \varphi (T^{-k}x_0 ) , v_n \big\rangle 
+ \Big\langle \sum_{j=0}^{k-1} \Psi^j \rho (T^{-(j+1)} x_0), v_n \Big\rangle\\
&=& \big\langle \Psi^k \varphi (T^{-k} x_0 ) , v_n \big\rangle 
+ \sum_{j=0}^{k-1} \rho_{n-j} (T^{-(j+1)} x_0).
\end{eqnarray*}
Since $\varphi$ is strongly continuous, the set $\Psi^k (\varphi(X))$ weakly converges 
to $\{ 0 \}$ in the Hausdorff sense. As a consequence, the first term above, namely  
$\langle \Psi^k  \varphi (T^{-k}x_0 ), v_n \rangle$, converges to zero 
as $n$ goes to infinite. Therefore, $\varphi$ has the following form:
\begin{equation}\label{la-formula}
\varphi(x) = \sum_{n \in \mathbb{Z}} \left( \sum_{j=0}^{\infty} \rho_{n-j} 
\big( T^{-(j+1)}x \big) \right) v_n.
\end{equation}

\vspace{0.1cm}

A closely related but slightly different case is that of a positive shift, that is, 
when $T: X \to X$ is a homeomorphism all of whose forward orbits are dense, and 
$\Psi$ is constant and coincides with the shift on the canonical basis $\{v_n\}$ 
of $\mathcal{H} \sim \ell^2(\mathbb{N}_0)$. (Notice that $\Psi$ is not 
surjective). Indeed, among all (non-neccessarily continuous) sections 
$\varphi \!: X \to \mathcal{H}$, there is a unique solution 
to (\ref{para-shift}), and its expression is given by 
\begin{equation}\label{bien}
\varphi(x) := 
\sum_{j=0}^{\infty} \left[\sum_{r=0}^{j} \rho_{j-r} \big( T^{-(r+1)} x \big)\right] v_j.
\end{equation}

To see that $\varphi(x)$ belongs to $\mathcal{H}$, we first 
claim that for all $y\in X$ and all $n \in \N$, we have 
\begin{equation}\label{des}
\big\| I(n,y)0 \big\| \leq 2C.
\end{equation}
Indeed, letting $M$ be the closure of the (forward) orbit of $ (x_0, 0)$ (where 
the topology on $\mathcal{H}$ is the weak one), we have $\|I(n,y) (v) \| \leq C$ for all 
$(y, v)\in M$ and all $n \in \N$. Since the forward orbits of $T$ are dense, each fiber 
$M_y$ is nonempty, hence we may take $v = v(y) \in M_y$. Using the triangular 
inequality and the fact that $I(n,y)$ is an isometry, we get 
$$\|I(n,y) 0 \| \leq \| v \|+\|I(n,y) v \| \leq 2C.$$
Now, a simple computation yields
$$I(n,y)0 = \sum_{r=0}^{n-1} \Psi^r \rho\big(T^{-(r+1)}(T^n y)\big)$$
for all $n \in \N$ and all $y\in X$. Using (\ref{des}), we obtain 
\begin{equation}\label{uno}
\sum_{j=0}^\infty\left|\sum_{r=0}^{n-1}\rho_{j-r}\big(T^{-(r+1)}(T^n y)\big)\right|^2 \leq 4C^2,
\end{equation}
where we let $\rho_{k} \equiv 0$ for $k < 0$. Given $N\in \N$, 
we choose $n > N$ and $y = T^{-n} x$ in (\ref{uno}), and we obtain 
$$\sum_{j=0}^N\left|\sum_{r=0}^{j}\rho_{j-r}\big(T^{-(r+1)} x \big)\right|^2 \leq 4C^2.$$
Since this holds for all $N \in \N$, we finally get $\|\varphi(x)\| \leq 2C$. In particular, 
$\varphi(x)$ belongs to $\mathcal{H}$.

To see that $\varphi$ is skew invariant, we just compute:
\begin{eqnarray*}
I(1,x) \varphi(x) 
&=& \Psi \varphi(x) + \rho(x)\\
&=& \Psi \left(\sum_{j=0}^{\infty}\left[\sum_{r=0}^{j}
\rho_{j-r} \big( T^{-(r+1)}x \big)\right]v_j\right)+\sum_{j=0}^{\infty}\rho_j(x) v_j\\
&=& \sum_{j=0}^{\infty}\left[\sum_{r=0}^{j}
\rho_{j-r} \big( T^{-(r+1)}x \big)\right]v_{j+1} + \sum_{j=0}^{\infty}\rho_j(x) v_j\\
&=& \sum_{j=1}^{\infty}\left[\sum_{r=0}^{j-1}
\rho_{j-1-r} \big( T^{-(r+1)}x \big) \right]v_{j} + \sum_{j=0}^{\infty}\rho_j(x) v_j\\
&=& \sum_{j=1}^{\infty}\left[\sum_{r=0}^{j-1}
\rho_{j-1-r} \big( T^{-(r+1)}x \big) \right]v_{j} + \sum_{j=0}^{\infty}\rho_j(x) v_j\\
&=& \sum_{j=1}^{\infty}\left[\sum_{r=1}^{j}
\rho_{j-r} \big( T^{-r}x \big) \right]v_{j} + \sum_{j=0}^{\infty}\rho_j(x) v_j\\
&=& \sum_{j=0}^{\infty}\left[\sum_{r=0}^j 
\rho_{j-r}\big(T^{-(r+1)} (Tx) \big)\right] v_j\\
&=& \varphi(T x ).
\end{eqnarray*}

To see that $\varphi$ is the unique skew-invariant function, 
we consider another such a function $\varphi^*: X \to \mathcal{H}$. 
For all $x \in X$ we have
$$\varphi^* ( T x ) - \varphi ( T x ) 
= I(1,x) \varphi^*(x) - I(1,x) \varphi(x)  
= \Psi \varphi^*(x) - \Psi \varphi(x),$$ 
that is, $(\varphi^* - \varphi)(Tx) = \Psi (\varphi^{*} - \varphi)(x)$. Defining 
$\phi_j \! : X \to \R$ by letting 
$$\sum_{j=0}^{\infty} {\phi_j}(x) v_j := \varphi(x) - \varphi^*(x),$$
this yields \esp 
$\phi_j (x) = \phi_{j-n} (T^{-n}x)$ \esp 
for all $n \in \N$. For $n > j$, this implies that $\phi_j(x) = 0$, 
hence $\varphi^* = \varphi$. 

Finally, since we know that there exists a continuous skew-invariant section, 
the expression (\ref{bien}) defines a continuous function. 

\vspace{0.02cm}

\begin{obs} {\em Since we know that the map $x \mapsto ctr(M_x)$ is skew invariant for 
any skew-invariant bounded subset $M \subset X \times \mathcal{H}$ whose projection 
on the first coordinate is onto, the vector $\varphi(x)$ above must coincide 
with $ctr (M_x)$ for all $x \in X$.} 
\end{obs} 



\vspace{0.3cm}

\noindent{\bf Acknowledgments.} A.~Navas would like to strongly thank Indira Chatterji 
for several discussions at the origin of Theorems C and D during a visit to the University of 
Orl\'eans (October 2010) for which he also acknowledges the invitation of the MAPMO Laboratory; 
moreover, he would like to thank Tsachik Gelander for his explanations on the Ryll-Nardewski 
theorem, \'Etienne Ghys for many remarks and comments, and Pierre Py for the suggestion 
of the Corollary to Theorem~C. M.~Ponce would like to thank Jean-Christophe Yoccoz and 
Raph\"ael Krikorian for their useful remarks. 

D.~Coronel was funded by the Fondecyt Post-doctoral Grant 3100092, 
A.~Navas was funded by the PBCT-Conicyt Research Project ADI-17 
and the Math-AMSUD Project DySET, and M.~Ponce was funded by 
the Fondecyt Grant 11090003 and the Math-AMSUD Project DySET.





\begin{small}

\end{small}


\begin{footnotesize}

\vspace{0.25cm}

\noindent{Daniel Coronel}

\noindent{Facultad de Matem\'aticas, PUC}

\noindent{Casilla 306, Santiago 22, Chile}

\noindent{E-mail: acoronel@mat.puc.cl}\\ 

\noindent{Andr\'es Navas}

\noindent{Dpto de Matem\'atica y C.C., USACH}

\noindent{Alameda 3363, Estaci\'on Central, Santiago, Chile}

\noindent{E-mail: andres.navas@usach.cl}\\

\noindent{Mario Ponce}

\noindent{Facultad de Matem\'aticas, PUC}

\noindent{Casilla 306, Santiago 22, Chile}

\noindent{E-mail: mponcea@mat.puc.cl}

\end{footnotesize}

\end{document}